[1994/12/01]

\documentclass[twoside,reqno,10pt]{amsart}
 \issueinfo{\textbf{2} (2002)}
 {}
 {}
 {2002}
 \pagespan{1}{39} \PII{}
 \copyrightinfo{2002}{D. Novikov, S. Yakovenko}

\usepackage{amsmath,amssymb,amsthm}
\usepackage[
    pdftitle={Quasialgebraicity of Picard--Vessiot extensions},
    pdfauthor={Dmitry Novikov, Sergei Yakovenko},
    bookmarksopen=false,
    pdfstartview=FitH,
    linkcolor=blue,
    citecolor=blue,
    urlcolor=blue,
     plainpages,
    ]{hyperref}

\def\arxivno#1{\texttt{#1}}

\def\arxivno#1{\href{http://arxiv.org/abs/#1}{\texttt{#1}}}
\def\MR#1{MR~\href{http://www.ams.org/mathscinet-getitem?mr=#1}{\textbf{#1}}}

\numberwithin{equation}{section}

\newtheorem{Cor}{Corollary}
\newtheorem{Thm}{Theorem}
\newtheorem{Prop}{Proposition}
\newtheorem{Lem}{Lemma}

\newtheorem{OtherThm}{Theorem}

\theoremstyle{definition}
\newtheorem{Def}{Definition}
\newtheorem{Ex}{Example}

\theoremstyle{remark}
\newtheorem{Rem}{Remark}

\renewcommand\ge{\geqslant}
\renewcommand\le{\leqslant}
\let\tildeaccent=\~
\renewcommand\~[1]{\widetilde{#1}}

\def\<{\left<}
\def\>{\right>}
\def\({\ifmmode\left(\else\textup{(}\fi}
\def\){\ifmmode\right)\else\textup{)}\fi}

\def\const{\operatorname{const}}

\def\Mat{\operatorname{Mat}}
\def\crit{\operatorname{crit}}

\def\Re{\operatorname{Re}}
\def\Im{\operatorname{Im}}
\def\Arg{\operatorname{Arg}}
\def\dist{\operatorname{dist}}

\def\diag{\operatorname{diag}}
\def\GL{\operatorname{GL}}

\renewcommand\:{\colon}
\def\R{{\mathbb R}}
\def\C{{\mathbb C}}

\def\e{\varepsilon}
\let\ssm=\smallsetminus
\def\pd#1#2{\frac{\partial#1}{\partial#2}}

\let\paragraph=\S
\def\secref#1{\paragraph\ref{#1}}

\def\S{\varSigma}
\def\CP{{\C P}}
\def\f{\varphi}
\def\N{{\mathbb N}}
\def\l{\lambda}

\begin{document}

\title
 { Quasialgebraicity of Picard--Vessiot fields }

\author
 {Dmitry Novikov \and Sergei Yakovenko}

\thanks{The research was supported by the Israeli Science Foundation grant
no.~18-00/1, the Killam grant of P.~Milman and the James S.
McDonnell Foundation.
\endgraf The paper is available online as
ArXiv preprint \href{http://arxiv.org/abs/math.DS/0203210}{{\tt
math.DS/0203210}}}

\address{Department of Mathematics\\
 University of Toronto\\
 Canada}

\curraddr{Department of Mathematics \\Purdue University \\
 West Lafayette \\ U.S.A.}

\email{{\tt dmitry@math.purdue.edu}}

\address{Department of  Mathematics\\
 The Weizmann Institute of Science\\
 Rehovot\\Israel}
\email{\href{mailto:yakov@wisdom.weizmann.ac.il}{\tt
yakov@wisdom.weizmann.ac.il}
\endgraf{\it WWW page\/}:
\href{http://www.wisdom.weizmann.ac.il/~yakov}{\tt
http://www.wisdom.weizmann.ac.il/\char'176 yakov/index.html}}
\date{June 2, 2002}


\dedicatory{To Vladimir Igorevich Arnold with admiration}

\begin{abstract}
We prove that under certain spectral assumptions on the monodromy
group, solutions of {Fuchsian systems} of linear equations on the
Riemann sphere admit explicit global bounds on the number of their
isolated zeros.
\end{abstract}

\maketitle

\tableofcontents

\section{Introduction}

\subsection{Fuchsian systems, monodromy}
Consider a system of first order linear ordinary differential
equations with rational coefficients on the complex Riemann sphere
$\CP^1$. In the matrix form such a system can be written as
\begin{equation}\label{ls}
  \dot X(t)=A(t)X(t),\qquad X\in\Mat_{n\times n}(\C),
\end{equation}
where $A(t)$ is the rational \emph{coefficients matrix} and $X(t)$
a multivalued \emph{fundamental matrix solution}, $\det
X(t)\not\equiv0$.

Assume that the system is \emph{Fuchsian}, that is, the
matrix-valued differential 1-form $\Omega=A(t)\,dt$ has only
simple poles on $\CP^1$, eventually including the point
$t=\infty$. Then the coefficients matrix $A(t)$ has the form
\begin{equation}\label{A}
  A(t)=\sum_{j=1}^m \frac{A_j}{t-t_j},\qquad A_j\in\Mat_{n\times n}(\C),
\end{equation}
with the constant \emph{matrix residues} $A_1,\dots,A_m$ at the
singular points $t_1,\dots,t_m$. The point $t=\infty$ is singular
if $A_1+\cdots+A_m\ne 0$, and then the residue at $t=\infty$ is
the negative of the above sum of finite residues.

\begin{Def}
The \emph{height} of the rational matrix function $A(t)$ as in
\eqref{A}, is the sum of norms of the residues,
\begin{equation}\label{res-norm}
  \|A(\cdot)\|=\|A_1\|+\cdots+\|A_m\|+\|A_\infty\|,
  \qquad A_\infty=-(A_1+\cdots+A_m).
\end{equation}
\end{Def}

The \emph{singular}, or \emph{polar locus} $\S=\{t_1,\dots,t_m\}$
of the system \eqref{ls} is the union of all \emph{finite}
singular points; if $t=\infty$ is a singular point also, then we
will denote $\S^*=\S\cup\{\infty\}\subset\CP^1$.

Any fundamental matrix solution $X(t)$ of the linear system
\eqref{ls} with a rational matrix of coefficients is an analytic
function on $\C$ ramified over the polar locus $\S$. If $\gamma$
is a closed loop avoiding $\S$, then analytic continuation
$\Delta_\gamma X(t)$ of the fundamental matrix solution $X(t)$
along $\gamma$ produces another fundamental solution, necessarily
of the form $X(t)M_\gamma$, where the constant invertible matrix
$M_\gamma\in\Mat_{n\times n}(\C)$ is called the \emph{monodromy
factor}. Choosing a different fundamental solution results in
replacing $M_\gamma$ by a matrix $CM_\gamma C^{-1}$ conjugate to
$M_\gamma$.

\begin{Def}\label{def:spectral-loop}
The linear system \eqref{ls} with a rational matrix $A(t)$ and the
singular locus $\S\subset\C$, is said to satisfy the
\emph{simple-loop spectral condition}, if all eigenvalues of each
monodromy factor $M_\gamma$ associated with any simple loop
(non-selfintersecting closed Jordan curve), belong to the unit
circle:
\begin{equation}\label{spec}
  \text{for any simple loop }\gamma\in\pi_1(\C\ssm\S),
  \quad\operatorname{Spec}M_\gamma\subset
  \{|\l|=1\}.
\end{equation}
\end{Def}

Clearly, this condition does not depend neither on the choice of a
fundamental solution, nor on the parametrization or even the
orientation of the loop.

\subsection{Principal result}
Let $T\subset\C\ssm\S$ be a simply connected domain in $\C\ssm\S$.
Then any fundamental matrix solution $X(t)$ is analytic in $T$
and, moreover, the linear space spanned by all entries
$x_{ij}(t)$, $i,j=1,\dots,n$, is independent of the choice of $X$.

The main problem addressed in the article is to place an explicit
upper bound on the number of isolated zeros of an arbitrary
function from this linear space.

One can rather easily see that this bound must necessarily depend
on the dimension $n$ of the system and the degree $m$ of the
rational matrix function $A(t)$. Simple examples suggest that the
height $r$ is also a relevant parameter. Besides, one should
impose certain restrictions on the spectra of the residue matrices
$A_i$ and exclude from consideration simply connected domains that
are ``too spiralling'' around one or more singularities (all these
examples are discussed in details in the lecture notes
\cite{montreal}).

It turns out that besides those already mentioned, there are no
other parameters that may affect the number of isolated zeros.
Moreover, an upper bound can be explicitly \emph{computed} in
terms of the parameters $n,m,r$. Recall that an integer valued
function of one or more integer arguments is \emph{primitive
recursive}, if it can be defined by several inductive rules, each
involving induction in only one variable. This is a strongest form
of computability, see \cite{manin,montreal}.

\begin{Thm}\label{thm:first}
There exists a primitive recursive function $\mathfrak N(n,m,r)$
of three natural arguments, with the following property.

If the Fuchsian system \eqref{ls}--\eqref{A} of dimension $n\times
n$ with $m$ finite singular points satisfies the simple-loop
spectral condition \eqref{spec} and its height is no greater than
$r$, then any linear combination of entries of a fundamental
solution in any triangular domain $T\subset\C\ssm\S$ has no more
than $\mathfrak N(n,m,r)$ isolated roots there.
\end{Thm}

The ``existence'' assertion concerning the primitive recursive
counting function $\mathfrak N$, is constructive. The proof of
Theorem~\ref{thm:first} in fact yields an algorithm for computing
$\mathfrak N(n,m,r)$ for any input $(n,m,r)\in\mathbb N^3$.

It is important to stress that the bound established in
Theorem~\ref{thm:first} is uniform over all triangles, all
possible configurations of $m$ distinct singular points, and all
combinations of residues $A_j$ of total norm $\le r$, provided
that the simple-loop spectral condition holds.

The assumption of triangularity of the domain $T$, as well as the
focus on linear combinations only, are not important. The general
assertion, Theorem~\ref{thm:main} on \emph{quasialgebraicity} of
Picard--Vessiot fields for Fuchsian systems meeting the spectral
condition \eqref{spec}, will be formulated in
\secref{sec:quasialg-def} after introducing all technical
definitions.

\subsection{The Euler system}\label{sec:euler}
The simplest example of a Fuchsian system is the \emph{Euler
system},
\begin{equation}\label{euler}
  \dot X=t^{-1}AX,\qquad A\in\Mat_{n\times n}(\C).
\end{equation}

The system \eqref{euler} can be explicitly solved:
$X(t)=t^A=\exp(A\ln t)$. The monodromy group of the Euler system
is cyclic and generated by the monodromy operator $\exp 2\pi iA$
for the simple loop encircling the origin. Taking $A$ in the
Jordan normal form (without losing generality) allows to compute
the matrix exponent and verify that the linear space spanned by
the components of $X(t)$, consists of \emph{quasipolynomials},
functions of the form
\begin{equation}\label{quasipol}
  f(t)=\sum_{\l\in\Lambda}t^\l\,p_\l(\ln t),
  \qquad\Lambda=\{\l_1,\dots,\l_n\}\subset\C,
  \quad p_\l\in\C[\ln t].
\end{equation}
The \emph{spectrum} $\Lambda$ of these quasipolynomials coincides
with the spectrum of the residue matrix $A$ and the  \emph{degree}
$d=\sum_\Lambda(1+\deg p_\l)$ is equal to the dimension $n$ of the
initial system.

Distribution of complex zeros of quasipolynomials depends on their
degree and spectrum. For example, if $\Lambda=\{\pm i\}$, then the
quasipolynomial $f(t)=t^i+t^{-i}=2\cos\ln t$ has an infinite
number of positive real roots accumulating to the origin.

On the other hand, if $\Lambda\subset\R$, then at least for
quasipolynomials with real coefficients, i.e., when $p_\l\in\R[\ln
t]$, the number of real positive zeros is no greater than $d-1$
(similarly to the usual polynomials corresponding to the case
$\Lambda\subset\N$).

It turns out that if $\Lambda\subset\R$ then not only real, but
also all complex roots can be counted. The following result, being
a particular case of Theorem~\ref{thm:first}, will be used as a
basis for the inductive proof of the general case as well.

\begin{Lem}[see \cite{jdcs-96a}, Theorem~2]\label{lem:euler-root}
The number of isolated roots of any quasipolynomial
\eqref{quasipol} of degree $\le d$ with arbitrary complex
coefficients but the \emph{real} spectrum
$\Lambda\subset[-r,r]\subset\R$, in any triangular domain $T$ not
containing the origin $t=0$, is explicitly bounded in terms of $d$
and $r$.

More precisely, in \cite{jdcs-96a} is shown that this number never
exceeds $4r+d-1$.\qed
\end{Lem}

\subsection{Systems with pairwise distant singular
points}\label{sec:distant-intro} The most difficult part of the
proof of Theorem~\ref{thm:first} is to treat \emph{confluent
singularities}, ensuring that the bounds on the number of zeros
would remain uniform even when the distances $|t_i-t_j|$ are
arbitrarily small. If this is not the case, more precisely, if the
bounds are allowed to depend on the configuration of the singular
points, then both the formulation and the proof can be
considerably simplified.

The corresponding result, formulated in a proper context below
(Theorem~\ref{thm:distant}, see \secref{sec:distant}), differs
from its unrestricted counterpart, Theorem~\ref{thm:second}, by
several instances.

\begin{enumerate}
    \item Fewer restrictions are imposed on the monodromy group:
    only \emph{small loops} around singular points must satisfy
    the spectral condition \eqref{spec}.
    \item Unlike the general case, the spectral condition for such
    small loops can be immediately verified by inspection of the
    eigenvalues of $A_j$ of the system: it is sufficient to
    require that all these eigenvalues must be real.
    \item Since the proof is considerably simplified, the bounds
    in this case are much less excessive, though still unlikely to
    be accurate.
    \item The price one has to pay is that the bounds for the number
    of zeros depend explicitly on the minimal distance between the
    singular points on $\C P^1$. These bounds explode as this
    distance tends to zero.
\end{enumerate}

This last observation on explosion of bounds notwithstanding, any
three points of $\C P^1$ can always be placed by an appropriate
conformal isomorphism to $0,1,\infty$. This implies computability
of the number of zeros for systems having only three Fuchsian
singular points on the sphere $\C P^1$, provided that the three
residue matrices have only real eigenvalues (cf.~with
\secref{sec:euler}). Such result may be of interest for the theory
of special functions, many of which are defined by this type of
systems.

\subsection{Applications to tangential Hilbert problem}
The problem on zeros of functions defined by Fuchsian equations,
is intimately related to the \emph{tangential}, or
\emph{infinitesimal Hilbert problem}
\cite{arnold-oleinik,arnold:quel-prob-sys-dyn}, see
\cite{montreal} for the detailed exposition and bibliography.
Recall that the problem concerns the maximal possible number of
isolated real zeros of Abelian integrals of the form
$I(t)=\oint_{H=t}\omega$, where $H$ is a bivariate polynomial and
$\omega$ a polynomial 1-form of a given degree $d$. The key
circumstance is the fact that Abelian integrals satisfy a more or
less explicitly known linear system of Picard--Fuchs differential
equations \cite{odo-2,redundant}.

The results formulated below, imply \emph{computability} of this
bound for polynomials $H$ whose critical points are sufficiently
distant from each other (Theorem~\ref{thm:distant}). The general
case still requires additional efforts to treat confluent
singularities, as explained in \cite{redundant}.

Yet there is one specific \emph{hyperelliptic} case when
$H(x,y)=y^2+F(x)$, $F\in \C[x]$. Computability of the bound in
this case under the additional assumption that all critical values
of $F$ are real, was proved in \cite{era-99} by almost the same
construction as exposed below. The difference concerns only two
technical instances. The preliminary folding of the Fuchsian
system (Appendix~\ref{sec:folding}) is made obsolete by the above
additional assumption. On the other hand, the isomonodromic
surgery explained in \secref{sec:surgery} and
Appendix~\ref{sec:factorization} in this case can be bypassed by a
suitably adapted application of Lyashko--Looijenga theorem
\cite{looijenga}.

\subsection*{Acknowledgements}
We are grateful to many people who explained us fine points of
numerous classical results from analysis, algebra and logic that
are used in the proof. Our sincere thanks go to
    V.~Arnol$'$d,
    J.~Bernstein,
    A.~Bolibruch,
    A.~Gabri\-\'elov,
    S.~Gusein-Zade,
    D.~Harel,
   Yu.~Ilyashenko,
    V.~Katsnelson,
    A.~Khovanski\u\i,
    V.~Matsaev,
    C.~Miller,
    P.~Milman,
    A.~Shen,
    Y.~Yomdin.

This work was partially done during the second author's stay at
the University of Toronto.

\section{Quasialgebraicity}

\subsection{Strategy}\label{sec:strategy}
The proof of Theorem~\ref{thm:first} is organized as follows.
First, we reformulate the assertion on zeros of linear
combinations of components of a fundamental solution, as a general
claim (Theorem~\ref{thm:second}) on \emph{quasialgebraicity} of
\emph{Picard--Vessiot fields}. Speaking loosely, this
quasialgebraicity means algorithmic computability of bounds for
the number of zeros of any function belonging to the function
field $\C(X)$ generated by the entries $x_{ij}(t)$ of any
fundamental solution $X(t)$.

The proof of quasialgebraicity goes essentially by induction in
the number of finite singular points of the system \eqref{ls}.
Actually, the assumptions on the system that are required (and
reproduce themselves in the induction) slightly differ from those
mentioned in Theorem~\ref{thm:first}. Instead of requiring all
singularities to be Fuchsian, we allow one of them (at infinity)
to be regular non-Fuchsian, with bounded negative Laurent matrix
coefficients. On the contrary, we require that all finite
singularities be aligned along the real axis and moreover belong
to the segment $[-1,1]$. The section \secref{sec:reduction}
contains the initial reduction, the proof of the fact that this
additional condition can be always achieved by suitable
transformations preserving the quasialgebraicity.

The hard core of the inductive step is provided by the following
general \emph{isomonodromic reduction principle}: if two given
Fuchsian systems have the same monodromy in a simply connected
polygonal domain $U\subset\C$ whose boundary is sufficiently
distant from all singularities, then the Picard--Vessiot fields
for these two systems, after restriction on $U$, are both
quasialgebraic or both not quasialgebraic simultaneously. The
accurate formulation of this principle is given in
\secref{sec:isomonodromic-reduction}.

The proof of the isomonodromic reduction rests upon another quite
general fact, the \emph{bounded meandering principle}
\cite{lleida,annalif-99}. It provides a possibility to majorize
explicitly the variation of argument of functions from
Picard--Vessiot fields, along segments distant from the singular
locus of a rational system. The accurate formulations are given in
\secref{sec:index} while the proofs (completely independent from
the rest of the paper) are moved to Appendix~\ref{sec:meandering}.

To carry out the inductive process, it is necessary to construct a
system with fewer singular points, isomonodromic with a given
system in a prescribed domain $U$. This \emph{isomonodromic
surgery} is closely related to the Hilbert 21st problem (the
\emph{Riemann--Hilbert problem}), see \cite{bolibr:kniga}. Its
analytic core is the quantitative version of the matrix
factorization problem, see Appendix~\ref{sec:factorization}. It is
this step that actually requires extending the class of Fuchsian
systems to that allowing one regular non-Fuchsian point (see
above).

\subsection{Picard--Vessiot fields: definitions, notations}
 \label{sec:quasialg-def}
Consider a linear system \eqref{ls} with rational coefficients
matrix $A(t)$ and let $X(t)$ be a fundamental matrix solution of
this system.

Consider the extension of the polynomial ring $\C[t]$ obtained by
adjoining all entries $x_{ij}(t)$ of the chosen fundamental
solution $X$: we will denote it by $\C[X]$ instead of a more
accurate but certainly more cumbersome notation
$\C[t,x_{11}(t),\dots,x_{nn}(t)]$. Moreover, from now on we will
assume (to simplify the language) that the list of generators
$X=\{x_{ij}\}$ of any Picard--Vessiot contains the independent
variable $t$.

The ring $\C[X]$ can be identified with a subring of the ring
$\mathcal O(t_*)$ of germs of analytic functions at any
nonsingular point $t_*\notin\S$. The field of fractions of $\C[X]$
will be denoted by $\C(X)$ and can be identified with a subfield
of the field $\mathcal M(t_*)$ of germs meromorphic at $t_*$. The
following properties of the ring $\C[X]$ and $\C(X)$ are obvious.

\begin{enumerate}
\item The construction of $\C(X)$ does not depend on the choice of the
fundamental solution $X$. For any two different nonsingular points
the ``realizations'' of $\C(X)$ by fields of meromorphic germs at
these points, are isomorphic to each other. However, this
isomorphism (e.g., by analytic continuation along a path
connecting these points, see below) is not canonical.

\item Any element $f\in\C[X]$ (respectively, from $\C(X)$), if
identified with a germ at $t_*\notin\S$, can be analytically
continued as a holomorphic (resp., meromorphic) function along any
path avoiding $\S$. If the path is a closed loop beginning and
ending at $t_*$, then the result of such continuation again
belongs to $\C[X]$ (resp., $\C(X)$). The corresponding
\emph{monodromy operator} is an automorphism of the ring (field),
extending the linear transformation $X\mapsto XM_\gamma$ on the
set of generators of $\C(X)$.

\item For any simply connected domain $U$ free from singular
points of the system \eqref{ls}, $\C[X]$ and $\C(X)$ can be
identified with a subring (resp., subfield) of the ring $\mathcal
O(U)$ of holomorphic (resp., the field $\mathcal M(U)$ of
meromorphic) functions in $U$. These subring and subfield will be
denoted respectively by $\C[X|U]$ and $\C(X|U)$ and referred to as
\emph{restrictions} of $\C[X]$ and $\C(X)$ on the domain $U$.

\item The field $\C(X)$ is closed by differentiation, and its
subfield of constants is $\C$. By definition this means that
$\C(X)$ is the \emph{Picard--Vessiot extension} of $\C(t)$, see
\cite{magid}. Somewhat surprisingly, we will not use this fact
directly.

\item The ring $\C[X]$ and the field $\C(X)$ are naturally
filtered by the degree. By definition, degree of an element
$f\in\C[X]$ is the lowest possible degree $d$ of any polynomial
expression $\sum c_{k,\alpha} t^k X^\alpha$, $c_{k,\alpha}\in\C$,
$|\alpha|+k\le d$, representing $f(t)$ (there may be several such
combinations, as we do not assume algebraic independence of the
entries $x_{ij}(t)$ of $X$). Degree of a fraction $f=g/h$,
$g,h\in\C[X]$, $f\in\C(X)$, is the maximum of $\deg g,\deg h$ for
the ``most economic'' representation of $f$.

\item The monodromy operators are degree-preserving.
\end{enumerate}

\subsection{Computability in the Picard--Vessiot fields}
To say about \emph{computability} of bounds on the number of
isolated zeros of functions from Picard--Vessiot fields, one has
in addition to the degree of functions introduce one or several
natural parameters characterizing the field itself (i.e., the
corresponding system \eqref{ls} for that matter), in terms of
which the bound should be expressed. Besides, when counting zeros
one should be aware of the multivaluedness of functions from
$\C(X)$.

The latter problem is resolved by the agreement to count zeros of
functions only in triangular domains free from singular points. To
address the former problem, we introduce  two admissible
\emph{standard classes} of systems \eqref{ls} with rational
coefficients and in each case we list explicitly the natural
parameters on which the bounds are allowed to depend.

\subsubsection{Fuchsian class}
Let $n,m,r$ be three natural numbers.

\begin{Def}
The \emph{Fuchsian class} $\mathcal F(n,m,r)$ consists of all
linear $n\times n$-systems \eqref{ls} with rational matrices
$A(t)$ of the form \eqref{A}, having no more than $m$ singular
points (including the one at infinity) and the height no greater
than $r$,
$$
\|A_1\|+\cdots+\|A_m\|+\|A_1+\cdots+A_m\|\le r.
$$
\end{Def}

This is the basic and the most important class. Its definition is
invariant by conformal changes of the independent variable
(M\"obius transformations). However, for technical reasons the
theorem on zeros has to be formulated and proved for
Picard--Vessiot fields built from another type of linear systems.

\subsubsection{The special class}
This class consists of systems having only one regular eventually
non-Fuchsian point at infinity. On the other hand, position of
finite Fuchsian singularities is subject to additional
restriction.

Let again $n,m,r$ be three natural numbers.

\begin{Def}\label{def:special}
The \emph{special} class $\mathcal S(n,m,r)$ consists of linear
systems \eqref{ls} with the matrix function $A(t)$ of the
following form,
\begin{equation}\label{special}
  A(t)=\sum_{j=1}^m\frac{A_j}{t-t_j}+\sum_{k=0}^r A_k'\,t^k,\qquad
  A_j,A_k'\in\Mat_{n\times n}(\C),
\end{equation}
such that:
\begin{enumerate}
  \item the finite singular points $t_1,\dots,t_m$ are all on the real
  interval $[-1,1]\subset\C$,
\begin{equation}\label{sc-poles}
  t_j\in\R,\quad |t_j|\le 1,\qquad\forall j=1,\dots,m,
\end{equation}
  \item all coefficients $A_j,A_k'$ are real matrices,
  $A_j,A_k'\in\Mat_{n\times n}(\R)$,
  \item the total norm of Laurent coefficients of
  the matrix $A(t)$ (including the singularity at infinity),
  is at most $r$,
\begin{equation}\label{sc-norm}
  \sum_{j=1}^m\|A_j\|+\sum_{k=0}^r\|A_k'\|\le r,
\end{equation}
  \item the singular point $t=\infty$ is regular, and
  any fundamental matrix solution $X(t)$
  grows no faster than $|t|^r$ as $t\to\infty$ along any ray:
\begin{equation}\label{sc-growth}
  \|X(t)\|+\|X^{-1}(t)\|\le c\,|t|^r,\qquad c>0,\ |t|\to+\infty,\ \Arg t=\const.
\end{equation}
\end{enumerate}
The parameter $r$ will be referred to as the \emph{height} of a
system from the special class.
\end{Def}

\begin{Rem}
The choice of parameters may seem artificial, in particular, the
assumption that $r$ simultaneously bounds so diverse things as the
residue norms, the growth exponent at infinity and the order of
pole of $A(t)$ at the regular singular point $t=\infty$. Yet since
our goal is only to prove computability without striving for
reasonable explicit formulas for the bounds, it is more convenient
to minimize the number of parameters, lumping together as many of
them as possible.
\end{Rem}

In both cases, the tuple of natural numbers describing the class,
will be referred to as \emph{the parameters} of the class or even
simply \emph{the parameters} of a particular system belonging to
that class. Accordingly, we will refer to Picard--Vessiot fields
built from solutions of the corresponding systems, as
Picard--Vessiot fields of Fuchsian or special type respectively,
the tuple $(n,m,r)$ being labelled as \emph{the parameters of the
field} (dimension, number of ramification points and the height
respectively). For the same reason we call $\S$ the ramification
locus of the Picard--Vessiot field $\C(X)$.

\subsubsection{Definition of quasialgebraicity}
Let $\C(X)$ be a Picard--Vessiot field of a Fuchsian or special
type. Recall that by $\C(X|T)$ we denote the restriction of the
field $\C(X)$ on a domain $T\subset\C P^1$.

\begin{Def}
The field $\C(X)$ is called \emph{quasialgebraic}, if the number
of isolated zeros and poles of any function $f\in\C(X|T)$ in any
\emph{closed} triangular domain $T\Subset\C\ssm\S$ free from
singular points of the corresponding system, is bounded by a
primitive recursive function $\mathfrak N'(d,n,m,r)$ of the degree
$d=\deg f$ and the parameters $n,m,r$ of the field.
\end{Def}

\begin{Rem}
It is important to stress that the bound for the number of zeros
is allowed to depend only on the parameters of the class as a
whole, and \emph{not} on specific choice of the system. In other
words, precisely as in Theorem~\ref{thm:first}, in any assertion
on quasialgebraicity of certain Picard--Vessiot fields the upper
bound for the number of isolated zeros and poles should be uniform
over:
\begin{enumerate}
  \item all configurations of $m$ singular points $t_1,\dots,t_m$,
  \item all triangles $T\subset\C\ssm\{t_1,\dots,t_m\}$,
  \item all matrix coefficients $A_j$, $A_k'$ meeting the
  restrictions on the total norm that occur in the definition of
  the standard classes,
  \item all functions $f$ of degree $\le d$ from $\C(X)$.
\end{enumerate}
\end{Rem}

\begin{Rem}
To avoid overstretching of the language, everywhere below outside
the principal formulations to be ``computable'' is synonymous to
be bounded by a primitive recursive function of the relevant
parameters (usually clear from the context). The triangles will be
always closed.
\end{Rem}

\subsection{Formulation of the results}
The principal result of this paper is the following apparent
generalization of Theorem~\ref{thm:first}.

\begin{Thm}[quasialgebraicity of Fuchsian fields]\label{thm:second}
A Picard--Vessiot field $\C(X)$ of Fuchsian type is quasialgebraic
if the corresponding  system \eqref{ls}--\eqref{A} satisfies the
spectral condition \eqref{spec}.
\end{Thm}

Actually this theorem is equivalent to Theorem~\ref{thm:first}. In
the case $m=1$ (i.e., for Euler systems) it follows from
Lemma~\ref{lem:euler-root}. In the general case
Theorem~\ref{thm:second} is obtained as a corollary to the similar
statement concerning the special class.

\begin{Thm}[main]\label{thm:main}
A Picard--Vessiot field $\C(X)$ of special type is quasialgebraic
if the corresponding  system \eqref{ls}, \eqref{special} satisfies
the spectral condition \eqref{spec}.
\end{Thm}

\subsection{Remarks}
Several obvious remarks should be immediately made in connection
with these definitions and formulations.

\subsubsection{Zeros versus poles}\label{sec:zeros-poles}
To prove quasialgebraicity, it is sufficient to verify that any
function $f$ from the \emph{polynomial ring} $\C[X]$ admits a
computable bound for the number of zeros only. Indeed, for a
rational function written as a ratio of two polynomials, the poles
may occur only at the roots of the denominator, since the
generators $X$ are always analytic in $T$.

\subsubsection{Independence on generators}\label{sec:indep-gener}
Quasialgebraicity of the Picard--Vessiot fields is essentially
independent of the choice of generators. Indeed, if
$Y=\{y_1,\dots,y_k\}\subset\C(X)$ is another collection of
multivalued functions such that $\C(Y)=\C(X)$ and the degrees
$\deg y_i$ are bounded by $k\in\mathbb N$, then any function
$f\in\C(Y)$ of degree $d$ in $Y$ is a rational combination of $X$
of degree $\le d'=2kd$. Computability in terms of $d$ or $d'$ is
obviously equivalent.

\begin{Ex}\label{ex:inverse}
The field $\C(X^{-1})$ generated by entries  of the inverse matrix
$X^{-1}(t)$, coincides with $\C(X)$. Indeed, $\det X\in\C[X]$ and
therefore $\C(X^{-1})\subseteq\C(X)$, and the role of $X$ and
$X^{-1}$ is symmetric. Since $X^{-1}$ is analytic outside $\S$
(has no poles), to verify quasialgebraicity of $\C(X)$ it is
sufficient to study zeros of functions from $\C[X^{-1}]$. Note
that in this case $\deg X^{-1}=n$, hence the two degrees $\deg_X$
and $\deg_{X^{-1}}$ on the same field $\C(X)=\C(X^{-1})$ are
easily re-computable in terms of each other.
\end{Ex}

\subsubsection{Triangular and polygonal domains}
The restriction by triangular domains is purely technical. As
follows from Proposition~\ref{prop:triangulation} below, one can
instead choose to count zeros in any semialgebraic domains in
$\C\simeq\R^2$ of known ``complexity'' (e.g., polygonal domains
with a known number of edges). However, in this case the bound on
zeros necessarily depends on the complexity of the domains. For
instance, polygonal domains with a large number of edges may wind
around singularities, spreading thus through many sheets of the
Riemann surface (see \cite{montreal}).

\begin{Prop}\label{prop:triangulation}
Let $\S$ be a point set of $m$ distinct points in $\R^2$ and
$D\Subset\R^2\ssm\S$ a compact simply connected semialgebraic
domain bounded by finitely many real algebraic arcs of the degrees
$d_1,d_2,\dots,d_k$ of total degree $d=d_1+\cdots+d_k$.

Then $D$ can be subdivided into at most $(d+\tfrac12)m(m+1)+2$
connected pieces such that each of them lies inside a triangular
domain free from points of $\S$.
\end{Prop}

\begin{proof}
Consider a partition of $\R^2\ssm\S$ into triangles, for instance,
slitting the plane $\R^2$ along some of the straight line segments
connecting points of $\S$. It can be also considered as a
spherical graph $G_1$ with $m+1$ vertices (one ``at infinity'')
and $M\le m(m+1)/2$ algebraic edges of degree $1$.

The domain $D$ together with its complement $\R^2\ssm D$ can be
also considered a spherical graph $G_2$ with $k$ vertices, $k$
edges and $2$ faces.

Superimposing these two graphs yields a new graph with new
vertices added at the points of intersection between the edges,
and the number of edges increased because the newly added vertices
subdivide some old edges. The number of such new vertices on each
edge, straight or ``curved'', can be immediately estimated from
B\'ezout theorem. As a result, the joint graph will have at most
$M(d_1+\cdots+d_k+1)=M(d+1)$ straight edges, at most
$(Md_i+1)+\cdots+(Md_k+1)=Md+k$ curvilinear algebraic edges and at
least $k$ vertices. By the Euler formula one can place an upper
bound on the number of faces of the joint graph, $F=2-k+E\le
2+M(2d+1)$.

Thus $D$ gets subdivided into a known number of connected
components, each of them belonging to exactly one triangle $T$ of
the triangulation of $\R^2\ssm\S$. This gives the required upper
bound on the total number of triangles covering $D$.
\end{proof}

Having proved this proposition, we can freely pass from one simply
connected domain to another, provided that they remain, say,
bounded by a known number of line segments or circular arcs.

\subsection{Relative quasialgebraicity}
In order to carry out the induction in the number of singular
points when proving Theorem~\ref{thm:main}, we need a relative
analog of quasialgebraicity of a Picard--Vessiot field in a given
domain $U$, eventually containing singularities. This technical
definition means a possibility of explicitly count zeros and poles
not ``everywhere in $\C$'' , but rather ``in $U$'' (with the
standard provision for multivalued functions to be restricted on
triangles). Unlike the global case, in the relative case the bound
is allowed to depend on the distance between the boundary of $U$
and the singular locus, yet this dependence must be computable.

This distance should take into account the eventual singular point
at infinity, so for any point set $S\subset\C$ we define
\begin{equation}\label{dist-infty}
  \dist(S,\S\cup\infty)=\inf_{t\in S,\ t_j\in
\S}\{|t-t_j|,~|t^{-1}|\}.
\end{equation}
If $S$ is compact and disjoint with $\S$, then this distance is
strictly positive.

Let $U\subset\C$ be a polygonal domain (for our purposes it would
be sufficient to consider only rectangles and their complements).
Denote by $\lceil x\rceil$ the integer part of a real number $x>0$
plus 1.

\begin{Def}
A Picard--Vessiot field $\C(X)$ built from solutions of a system
\eqref{ls} of Fuchsian or special type, is said to be
\emph{quasialgebraic in $U$}, if:
\begin{enumerate}
    \item the boundary $\partial U$ of the
    domain $U$ is bounded and contains no singularities of the
    system \eqref{ls}, so that
    $\dist(\partial U,\S\cup\infty)=\rho>0$,

    \item the number of isolated zeros and poles of any function
    $f\in\C(X)$ in any triangle $T\Subset U\ssm \S$ free from
    singular points, is bounded by a primitive recursive function
    $\mathfrak N''(d,n,m,r,s)$ of the degree $d=\deg f$, the
    parameters $n,m,r$ of the field and the inverse distance
    $s=\lceil 1/\rho\rceil\in\mathbb N$.
\end{enumerate}
\end{Def}

If the Riemann sphere $\C P^1$ is tiled by finitely many domains
$U_1,\dots,U_k$ and the field $\C(X)$ is quasialgebraic in each
$U_i$, then it is globally quasialgebraic provided that the
boundaries of the domains never pass too close to the
singularities.

\section{Alignment of singular points: the first reduction}
 \label{sec:reduction}

\subsection{Transformations of Picard--Vessiot fields}
In this section we show how Theorem~\ref{thm:main} implies
Theorems~\ref{thm:first} and~\ref{thm:second}.

This reduction consists in construction of an algebraic (ramified,
multivalued) change of the independent variable $z=\f(t)$ using
some rational function $\f\in\C(t)$, that ``transforms'' an
arbitrary Fuchsian system \eqref{ls}--\eqref{A} to another
Fuchsian system with respect to the new variable $z$, having all
singular points $z_j$ on $[-1,1]$ and only real residue matrices.
Speaking loosely, one has to ``fold'' the Riemann sphere $\C P^1$
until all singular points (both the initial singularities and the
singularities created when folding) fall on the unit segment.
Meanwhile the norms of residues of the obtained system should
remain explicitly bounded. The general idea of using such
transformation belongs to A.~Khovanski\u\i~\cite{asik:unpub}.

The word ``transforms'' appears in the quotation marks since the
ramified algebraic changes of the independent variable do not
preserve rationality of the coefficients. Indeed, if $X(t)$ is a
(germ of) fundamental matrix solution to a linear system
\eqref{ls} with rational matrix of coefficients $A(t)$ and
$z=\f(t)$ is a rational map of the $t$-sphere onto the $z$-sphere,
the transform\footnote{The prime $'$ never in this article is used
to denote the derivative.} $X'(z)=X(\f^{-1}(z))$ of $X(t)$ by the
\emph{algebraic} (\emph{non}-rational) change of the independent
variable $\psi=\f^{-1}$ does not in general satisfy a system of
equations with rational coefficients. Performing the change of
variables, we obtain
\begin{equation*}
  \frac{dX'(z)}{dz}=
 \(1\bigg/{\frac{d\f}{dz}(\f^{-1}(z))}\)A(\f^{-1}(z))X'(z),
\end{equation*}
and see that the new matrix of coefficients has in general only
algebraic entries. This means that the field $\C(\psi^*X)$
generated by entries of the pullback $\psi^*X=X\circ\psi=X\circ
\f^{-1}$, is not a differential extension of $\C(z)$: it has to be
completed.

\begin{Rem}
For any multivalued function $f\in\C(X)$ ramified over $\S$, its
pullback (transform) $\psi^*f(z)=f(\f^{-1}(z))$ by $\psi=\f^{-1}$
is a multivalued function ramified over the union
$\S'=\f(\S)\cup\crit\f$ consisting of the \emph{direct} $\f$-image
of $\S$ and the critical locus $\crit\f$ over which $\f^{-1}$ is
ramified.
\end{Rem}

\begin{Lem}[Alignment of singularities of Fuchsian systems]\label{lem:folding}
For any Picard--Vessiot field $\C(X)$ of Fuchsian type, there can
be constructed an algebraic transformation $\psi=\f^{-1}$ inverse
to a rational map $\f\in\C(t)$, and a Picard--Vessiot field
$\C(Y)$, also of Fuchsian type, such that\textup:
\begin{enumerate}
    \item $\psi^*\C(X)=\C(\psi^* X)\subseteq\C(Y)$,
    \item for any $f\in\C(X)$ of degree $d$, the pullback
    $\psi^*f\in\C(Y)$ has degree $\le d$ relative to $\C(Y)$,
    \item the ramification locus of $\C(Y)$ consists of only finite real
    points, and the corresponding residue matrices are all real,
    \item the degree of the map $\f$ and the parameters of the field
    $\C(Y)$, in particular its height, are bounded by computable
    functions of the parameters of the field $\C(X)$.
\end{enumerate}
In addition, if the spectral condition \eqref{spec} held for the
initial Fuchsian system generating the field $\C(X)$, then it will
also hold for the Fuchsian system generating the field $\C(Y)$.
\end{Lem}

\subsection{Special vs.~Fuchsian systems: derivation of
 Theorem~\ref{thm:second} from Theorem~\ref{thm:main}}
Obviously, quasialgebraicity of $\C(Y)$, when proved, implies that
its subfield $\C(\psi^*X)$ is also quasialgebraic. By
Proposition~\ref{prop:triangulation}, this means also
quasialgebraicity of $\C(X)$. Indeed, $\f$-preimages and
$\f$-images of triangles are semialgebraic domains of explicitly
bounded complexity, provided the degree of the map $\f$ is known.

As soon as all singular points of a Fuchsian system are finite and
aligned along the real axis $\R$, by a suitable affine
transformation of the independent variable they can be brought to
the segment $[-1,1]$. This transformation does not affect the
residues of the Fuchsian system, therefore all additional
conditions from the definition of the special class can be
satisfied: in this case $A'_k=0$ and the growth exponent at
infinity is also zero.

The folding construction shows that the question on
quasialgebraicity of an arbitrary Picard--Vessiot field $\C(X)$ of
Fuchsian type can be reduced using Lemma~\ref{lem:folding} to the
question on quasialgebraicity of a certain auxiliary field of
special type, with explicitly bounded parameters.

The proof of Lemma~\ref{lem:folding} is completely independent of
the rest of the paper and is moved to Appendix~\ref{sec:folding}.
{From} this moment on we concentrate exclusively on the proof of
Theorem~\ref{thm:main}.

\section{Computability of index}\label{sec:index}

\subsection{Index}
For a multivalued analytic function $f(t)\not\equiv0$ ramified
over a finite locus $\S$, and an arbitrary piecewise-smooth
oriented path $\gamma\subset\C\ssm\S$ avoiding this locus, denote
by $V_\gamma(f)$ the index of $f$ along $\gamma$.

More precisely, if $f$ has no zeros on $\gamma$, then $\Arg f$
admits selection of a continuous branch along $\gamma$ and the
real number $V_\gamma(f)$ is defined as $\Arg f(\text{end})-\Arg
f(\text{start})$ (this number is an integer multiple of $2\pi$ if
$\gamma$ is closed). If $f$ has zeros on $\gamma$, then
$V_\gamma(f)$ is defined as
\begin{equation*}
  V_\gamma(f)=\limsup_{c\to 0}V_\gamma(f-c),
\end{equation*}
where the upper limit is taken over values $c\notin f(\gamma)$. By
definition, index of the identically zero function is set equal to
$0$ along any path.

\subsection{Index of segments distant from the critical locus}
Consider a linear system \eqref{ls}, \eqref{special} from the
special class $\mathcal S(n,m,r)$ with the singular locus $\S$ and
the corresponding Picard--Vessiot field $\C(X)$.

If $\gamma$ is a bounded rectilinear segment not passing through
singular points, then the distance
$\rho=\dist(\gamma,\S\cup\infty)$ from $\gamma$ to $\S\cup\infty$,
defined as in \eqref{dist-infty}, is positive.

\begin{Lem}\label{lem:ind-seg}
For any Picard--Vessiot field $\C(X)$ of special type the absolute
value of index of any function $f\in\C(X)$ along any rectilinear
segment $\gamma\subset\C\ssm\S$ is bounded by a computable
function of $d=\deg f$, the parameters of the field and the
inverse distance $s=\lceil 1/\rho\rceil$.
\end{Lem}

This result is a manifestation of the general \emph{bounded
meandering principle}, see \cite{lleida,annalif-99}. It is derived
in Appendix~\ref{sec:meandering} from theorems appearing in these
references.

\subsection{Index of small arcs}
As the segment $\gamma$ approaches the singular locus, the bound
for the index along $\gamma$ given by Lemma~\ref{lem:ind-seg}
explodes. However, for Picard--Vessiot fields of Fuchsian and
special types, the index along small circular arcs around all
finite singularities remains computable.

\begin{Lem}\label{lem:ind-circ}
Assume that the Picard--Vessiot field $\C(X)$ is of Fuchsian or
special type and $t_j\in\S$ is a ramification point.

Then there exist arbitrarily small circles around $t_j$ \(not
necessarily centered at $t_j$ but containing it strictly inside\)
with the following property. The absolute value of the index of
any function $f\in\C(X)$ of degree $d$ along any part of any such
circle is bounded by a computable function of $d$ and the
parameters of the field.
\end{Lem}

In other words, index along some (certainly not all!) small
circular arcs around every finite singularity is bounded by a
computable function of the available data. The proof is given in
Appendix~\ref{sec:meandering}.

\subsection{Corollaries on quasialgebraicity}\label{sec:quasialg-holom}
The above results on computability of the index immediately imply
some results on relative quasialgebraicity.

\begin{Cor}
Any Picard--Vessiot fields of Fuchsian or special type is
relatively quasialgebraic in any domain $U$ free from ramification
points.
\end{Cor}

\begin{proof}
Indeed, any ``polynomial'' $f\in\C[X]$ has an explicitly bounded
index along any triangle $T\subset U$ by Lemma~\ref{lem:ind-seg}.
By the argument principle, this index bounds the number of zeros.
The answer is given in terms of the inverse distance
$\lceil\dist^{-1}(\partial U,\S\cup\infty)\rceil$ and the
parameters of the field. The assertion on zeros and poles of
rational functions from $\C(X)$ follows from the remark in
\secref{sec:zeros-poles}.
\end{proof}

Another application concerns quasialgebraicity of
\emph{subfields}. Assume that $U$ is a polygon (say, triangle) and
$H=\{h_1,\dots,h_k\}\subset\C(X|U)\cap\mathcal M(U)$ is a
collection (list) of functions from $\C(X)$ that after restriction
on $U$ are single-valued (i.e., meromorphic) there. Without loss
of generality we may assume that $h_i$ are in fact holomorphic in
$U\ssm\S$, multiplying them if necessary by the appropriate
polynomials from $\C[t]\subset\C(X)$.

\begin{Cor}\label{cor:quasialg-merom}
If $\C(X)$ is a Picard--Vessiot field of Fuchsian or special type
and $H\subset\C(X)\cap\mathcal O(U)$ is a collection of branches
holomorphic in a triangular domain $U$, then the subfield
$\C(H)\subseteq\C(X)$ is quasialgebraic in $U$.
\end{Cor}

\begin{proof}
The same argument principle can applied to polynomials $f\in\C[H]$
and the domain obtained by deleting from $U$ small disks around
singularities, as described in Lemma~\ref{lem:ind-circ}. Since $f$
has no poles in such domain, Lemmas~\ref{lem:ind-seg}
and~\ref{lem:ind-circ} together yield a computable upper bound on
the number of zeros $f$ in any triangle $T\subset U$, even if the
sides of the triangle pass arbitrarily close to the singularities
inside $U$. The case of rational functions is treated as above.
\end{proof}

\begin{Rem}
The domain $U$ may be not a triangle but a polygon with any
apriori bounded number of sides. As usual, if the bound for
functions $f\in\C(H)$ is to be expressed in terms of the degree of
$f$ with respect to the initial generators, then the degree of $H$
(the maximum of degrees of all $h_i$) must be explicitly
specified, as explained in \secref{sec:indep-gener}.
\end{Rem}

\section{Isomonodromic reduction}\label{sec:isomonodromic-reduction}
The two corollaries from \secref{sec:quasialg-holom} suggest that
an obstruction to the (relative) quasialgebraicity lies in the
non-trivial monodromy of the Picard--Vessiot field. The results of
this section show that this is the \emph{only} obstruction: at
least for Picard--Vessiot fields of special type, the monodromy
group is the only factor determining whether the field is
quasialgebraic or not.

\subsection{Isomonodromic systems, isomonodromic fields}
Let $U\subseteq\C$ be a (polygonal) bounded domain and
\begin{equation}\label{2sys}
  \dot X=A(t)X,\quad \dot Y=B(t)Y,
  \qquad A,B\in\Mat_n(\C(t)),
\end{equation}
two linear systems of the same size with rational coefficients.
Denote by $\S_A$ and $\S_B$ their respective singular loci.

\begin{Def}
The two systems are said to be isomonodromic in the domain $U$,
if:
\begin{enumerate}
    \item they have the same singular points inside $U$, $\S_A\cap
    U=\S_B\cap U=\S$, and no singular points on the boundary
    $\partial U$, and
    \item all monodromy operators corresponding to loops
    \emph{entirely belonging to $U$}, are simultaneously
    conjugated for these two systems.
\end{enumerate}
Two Picard--Vessiot fields are said to be isomonodromic in $U$, if
the respective linear systems are isomonodromic there.
\end{Def}

Choosing appropriate fundamental solutions (e.g., their germs at a
nonsingular point $t_*\in U\ssm\S$), one can assume without loss
of generality that the monodromy factors are all equal:
\begin{equation*}
 \forall\gamma\subset U, \quad
 X^{-1}\cdot \Delta_\gamma
 X=Y^{-1}\cdot \Delta_\gamma Y.
\end{equation*}
The matrix ratio $H=XY^{-1}$ of \emph{these} solutions is
single-valued in $U$ (note that this may not be the case for
another pair of fundamental solutions). Note also that after
continuing $H$ along a path leaving $U$, we arrive at a different
branch of $H$ which may well be ramified even in $U$.

Let $U\subset\C$ be a rectangle with sides parallel to the real
and imaginary axes and having no singularities on the boundary.

\begin{Lem}\label{lem:isomon}
Two Picard--Vessiot fields $\C(X)$ and $\C(Y)$ of the special
type, isomonodromic in a rectangle $U$, are either both
quasialgebraic, or both not quasialgebraic in $U$.

More precisely, if one of two isomonodromic fields admits a
computable upper bound $\mathfrak N_1(d,n,m,r)$ for the number of
zeros \(poles\) in $U$, then the other also admits the similar
bound $\mathfrak N_2(d,n,m,r,s)$, additionally depending \(also in
a computable way\) on the inverse distance $s$ between $\partial
U$ and $\S_A\cup\S_B\cup\infty$.
\end{Lem}

The proof of the Lemma occupies the rest of this section. The
strategy is, assuming that $\C(X)$ is quasialgebraic, to prove
quasialgebraicity of the \emph{compositum} $\C(X,Y)$, the
Picard--Vessiot field generated by entries of both fundamental
solutions. Clearly, this is sufficient since
$\C(Y)\subseteq\C(X,Y)$ will in this case also be quasialgebraic.
To count zeros of a function from $\C(X,Y)$ we use a modification
of the popular derivation-division algorithm
\cite{umn-91,roussarie:h16} based on the Rolle lemma for real
functions. The role of the derivation in this modified algorithm
is played by the operator $\Im$ of taking the imaginary part.

\subsection{Real closeness}\label{sec:re-im}
First we observe that the Picard--Vessiot field of special type is
closed by taking real or imaginary parts. We say that the system
\eqref{ls} is \emph{real}, if the coefficients matrix $A(t)$ of
this system is real on $\R$.

\begin{Prop}\label{prop:reim}
Let $\sigma\subset\R\ssm\S$ be a real segment free from
singularities of a real system \eqref{ls}. Then for any function
$f\in\C(X)$ there exist two functions from $\C(X)$, denoted
respectively by $\Re_\sigma f$ and $\Im_\sigma f$, such that
\begin{equation*}
  \Re(f|_\sigma)=(\Re_\sigma f)|_\sigma, \qquad
  \Im(f|_\sigma)=(\Im_\sigma f)|_\sigma.
\end{equation*}
The degree of both functions does not exceed $\deg f$.
\end{Prop}

In other words, the real/imaginary part of any branch of any
function $f$ on any \emph{real} segment $\sigma$ extends again on
the whole set $\C\ssm\S$ as an analytic function from the same
field $\C(X)$ and the degree is not increased.

\begin{proof}
Since the assertion concerns only the field, it does not depend on
the choice of the fundamental solution $X(t)$. Choose this
solution subject to an initial condition $X(t_*)=E\in\GL(n,\R)$,
where $t_*\in\sigma$ is arbitrary point on the segment. Then the
solution remains real-valued on $\sigma$, $X(t)\in\GL(n,\R)$ for
any $t\in\sigma$, since the system \eqref{ls} is real. The trivial
operators $\Re_\sigma$, $\Im_\sigma$, defined on the generator set
$X$ of $\C(X)$ as
\begin{equation*}
  \Re_\sigma X=X,\quad \Im_\sigma X=0;\qquad \Re_\sigma c=\Re
  c,\quad \Im_\sigma c=\Im c\quad
  \forall c\in\C,
\end{equation*}
can be naturally and uniquely extended on the field $\C(X)$ of
rational combinations with constant complex coefficients. They
possess all the required properties.
\end{proof}

\begin{Rem}\label{rem:globalRe}
Note that, unlike single-valued functions, the result in general
depends on the choice of $\sigma$ (this explains the notation).
Yet if the restriction $h\in\C(X|U)$ is a single-valued branch in
a domain $U$ intersecting the real axis, then $\Re_\sigma h$ and
$\Im_\sigma h$ do not depend for the choice of a segment
$\sigma\subset U\cap \R$. This immediately follows from the
uniqueness theorem for meromorphic functions.
\end{Rem}

\subsection{Index and real part: the Petrov lemma}
Let $f\:\sigma\to\C$ be a complex-valued function without zeros on
a line segment $\sigma\subset\C$, and $g\:\sigma\to \R$ is the
imaginary part of $f$, a real-valued function on $\sigma$. We
consider the case when $g$ is real analytic on $\sigma$ and hence
has only isolated zeros. The following elementary but quite
powerful trick was conceived by G.~Petrov \cite{petrov:trick}.

\begin{Prop}\label{prop:petrov}
If $g=\Im f\not\equiv0$ has $k$ isolated zeros on $\sigma$, then
the index of $f$ on $\sigma$ is bounded in terms of $k$\textup:
\begin{equation*}
  |V_\sigma(f)|\le \pi(k+1).
\end{equation*}
\end{Prop}

\begin{proof}
This is an ``intermediate value theorem'' for the map of the
interval $\sigma$ to the circle, sending $t\in\sigma$ to
$f(t)/|f(t)|$. If $V_\sigma(f)\ge\pi$, then the image of $f$
should cover two antipodal points on the circle, and hence contain
one of the two arcs connecting them on the sphere. But any such
arc intersects both real and imaginary axis.

The general case is obtained by subdividing $\sigma$ on parts
along which variation of argument is equal to $\pi$.
\end{proof}

\subsection{Compositum of Picard--Vessiot fields}
If $\C(X)$ and $\C(Y)$ are two Picard--Vessiot fields generated by
fundamental solutions of the two systems \eqref{2sys}, their
compositum $\C(X,Y)$, the field generated by entries of both
fundamental solutions, is again a Picard--Vessiot field. This is
obvious: the block diagonal matrix $\diag\{X,Y\}$ satisfies a
linear system with the block diagonal coefficients matrix
$\diag\{A,B\}$. Notice that this system will be of Fuchsian
(resp., special) class if both systems \eqref{2sys} were from this
class. The dimension, number of singularities and the height of
the system generating the compositum field, are computable in
terms of the respective parameters of the original fields.

As a consequence, for any function $f\in\C(X,Y)$ its index along
segments or arcs described in Lemmas~\ref{lem:ind-seg}
and~\ref{lem:ind-circ}, is explicitly computable in terms of the
degree of $f$ and the parameters of the two fields.

\subsection{Demonstration of Lemma~\ref{lem:isomon}}
As was already mentioned, the goal would be achieved if we prove
that the compositum $\C(X,Y)$ in $U$ is quasialgebraic provided
that $\C(X)$ is quasialgebraic there. The proof begins with a
series of preliminary simplifications.

1. We assume that the fundamental solutions $X,Y$ are chosen with
coinciding monodromy factors, so that their matrix ratio
$H=YX^{-1}$ is single-valued hence meromorphic in $U$. Its degree
in $\C(X,Y)$ is $n+1$ (Example~\ref{ex:inverse}). Thus to prove
the Lemma, it is sufficient to prove quasialgebraicity of the
field $\C(X,H)$:
\begin{equation*}
  \C(Y)\subset\C(X,Y)=\C(X,YX^{-1})=\C(X,H).
\end{equation*}
Since the matrix function $H$ has no poles in $U\ssm\S$, to prove
quasialgebraicity of $\C(X,H)=\C(X)(H)$, it is sufficient to
majorize the number of zeros of functions from $\C(X)[H]$, i.e.,
polynomially depending on the entries of $H$ (cf.~with Remark in
\secref{sec:zeros-poles}).

2. By Remark~\ref{rem:globalRe}, $\Re H$ and $\Im H$ are
unambiguously defined matrix functions, holomorphic in $U\ssm\S$
and real on $U\cap\R$. Replacing $H$ by $\diag\{\Re H,\Im H\}$ if
necessary, we may assume that $H$ itself is real on $\R\cap U$.

3. Since $H$ has no poles in $U\ssm\S$, to prove quasialgebraicity
it is sufficient to consider only ``polynomials in $H$'', i.e.,
the ring $\C(X)[H]\subset\C(X)(H)=\C(X,H)$. By the assumed
quasialgebraicity of $\C(X)$, the number of poles of each such
function, represented as $\sum x_\alpha(t)H^\alpha$,
$x_\alpha\in\C(X)$, is computable in terms of $\deg
x_\alpha\le\deg f$, the parameters of the field $\C(X)$ and the
relative geometry of $U$ and $\S_A\cup\S_B$. Therefore it remains
only to majorize the number of zeros of functions from this ring.

4. Any function from $f=\C(X)[H]$ can be represented as the finite
sum
\begin{equation}\label{sum}
  f=\sum_{i=1}^\nu x_i h_i,\qquad x_i\in\C(X),\ h_i\in\C[H],\ \Im h_i|_\R=0.
\end{equation}
The number $\nu$ and the degrees $\deg_X x_i$, $\deg_H h_i$ are
all explicitly computable in terms of $\deg_{X,Y} f$ and $n$.

After all this preparatory work we can prove by induction in
$\nu$, the number of terms in \eqref{sum}, that the number of
isolated zeros of any function of such form in any triangle
$T\subset U\ssm\S$, is computable in terms of the parameters of
the fields $\C(X)$ and $\C(Y)$, the degrees $\deg_X x_i$ and
$\deg_H h_i$ and the distance from $U$ to
$\S_A\cup\S_B\cup\infty$. The proof is very much similar to the
popular ``differentiation--division algorithm'' used to majorize
the number of real zeros of fewnomials in one variable
\cite{umn-91,roussarie:h16}.

5. For $\nu=1$ the assertion is immediate. Zeros of $x_1 h_1$ in
any triangle $T\subset U\ssm\S$ are either zeros of $x_1\in\C(X)$
or zeros of $h_1$. The number of zeros of the first kind in any
triangle is bounded by virtue of quasialgebraicity of $\C(X)$. The
number of zeros of second kind is bounded by
Corollary~\ref{cor:quasialg-merom}, since $h_1$ is holomorphic
with computable index along the boundary. Thus the base of
induction is established.

6. For a sum \eqref{sum} involving $\nu>1$ terms, we first divide
$f$ by $x_\nu$. This division may result in a loss or acquisition
of a computable number of zeros or poles (occurring as poles and
zeros of $x_\nu$ respectively). The degree of $x_j$, $j\ne 1$, as
well as the degree of $f$ may increase at most by $\deg x_i\le d$.
After such division we may without loss of generality assume that
$x_\nu\equiv1$, that is,
\begin{equation*}
  f=h_\nu+\sum_{i=1}^{\nu-1} x_i h_i, \qquad h_1,\dots,h_\nu\in\C[H],
  \quad x_1,\dots,x_{\nu-1}\in\C(X).
\end{equation*}

%
%

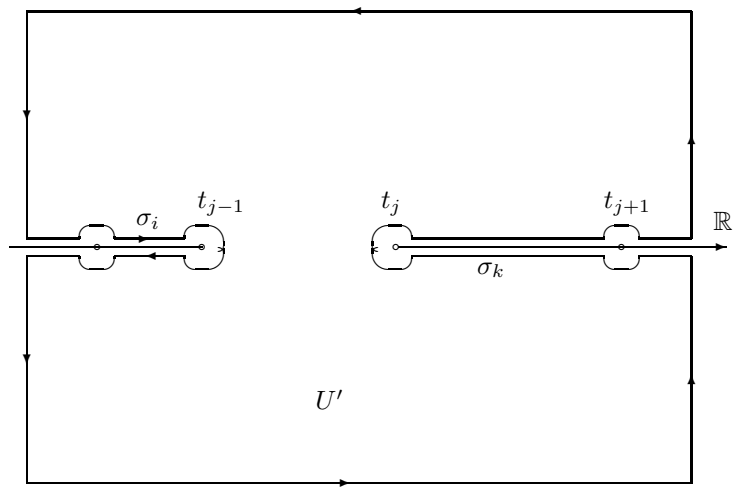
\begin{figure}
\setlength{\unitlength}{0.00030000in}%
\begin{picture}(12324,8124)(289,-7573)
\thinlines \put(1801,-3361){\oval(600,450)[tr]}
\put(1801,-3361){\oval(600,450)[tl]}
\put(1801,-3661){\oval(600,450)[bl]}
\put(1801,-3661){\oval(600,450)[br]}
\put(10801,-3361){\oval(600,450)[tr]}
\put(10801,-3361){\oval(600,450)[tl]}
\put(10801,-3661){\oval(600,450)[bl]}
\put(10801,-3661){\oval(600,450)[br]} \put(3976,-3484){\oval(  2,
54)[br]} \put(3628,-3484){\oval(698,698)[tr]}
\put(3628,-3361){\oval(654,452)[tl]}
\put(3628,-3661){\oval(654,452)[bl]}
\put(3628,-3538){\oval(698,698)[br]} \put(3976,-3538){\oval(  2,
54)[tr]} \put(6874,-3361){\oval(654,452)[tr]}
\put(6874,-3484){\oval(698,698)[tl]} \put(6526,-3484){\oval(  2,
54)[bl]} \put(6526,-3538){\oval(  2, 54)[tl]}
\put(6874,-3538){\oval(698,698)[bl]}
\put(6874,-3661){\oval(654,452)[br]} \put(1801,-3511){\circle{80}}
\put(3601,-3511){\circle{80}} \put(6926,-3511){\circle{80}}
\put(10801,-3511){\circle{80}} \put(601,-3361){\line( 1, 0){900}}
\put(2101,-3361){\line( 1, 0){1200}} \put(601,-3361){\line( 1,
0){900}} \put(601,-3661){\line( 1, 0){900}}
\put(2101,-3661){\line( 1, 0){1200}} \put(7201,-3361){\line( 1,
0){3300}} \put(7201,-3661){\line( 1, 0){3300}}
\put(11101,-3361){\line( 1, 0){900}} \put(11101,-3661){\line( 1,
0){900}} \put(601,-3361){\line( 0, 1){3900}} \put(601,539){\line(
1, 0){11400}} \put(12001,539){\line( 0,-1){3900}}
\put(601,-3661){\line( 0,-1){3900}} \put(601,-7561){\line( 1,
0){11400}} \put(12001,-7561){\line( 0, 1){3900}}
\put(6001,-7561){\vector( 1, 0){150}} \put(12001,-5761){\vector(
0, 1){ 75}} \put(12001,-1636){\vector( 0, 1){ 75}}
\put(6226,539){\vector(-1, 0){ 75}} \put(601,-1261){\vector(
0,-1){ 75}} \put(601,-5461){\vector( 0,-1){ 75}}
\put(2626,-3361){\vector( 1, 0){ 75}} \put(2701,-3661){\vector(-1,
0){ 75}} \put(301,-3511){\line( 1, 0){3300}}
\put(6976,-3511){\vector( 1, 0){5625}}
\put(5551,-6286){\makebox(0,0)[lb]{\smash{$U'$}}}
\put(12376,-3211){\makebox(0,0)[lb]{\smash{$\R$}}}
\put(2476,-3136){\makebox(0,0)[lb]{\smash{$\sigma_i$}}}
\put(8326,-3961){\makebox(0,0)[lb]{\smash{$\sigma_k$}}}
\put(6676,-2836){\makebox(0,0)[lb]{\smash{$t_j$}}}
\put(3526,-2836){\makebox(0,0)[lb]{\smash{$t_{j-1}$}}}
\put(10501,-2836){\makebox(0,0)[lb]{\smash{$t_{j+1}$}}}
\end{picture}
\caption{Slit rectangle $U$}
\end{figure}

Being convex, an arbitrary triangle $T\subset U\ssm\S$ may
intersect $\R$ by at most one connected segment: $T\cap\R\subseteq
\sigma\subset\R\ssm\S$, where $\sigma$ is one of the intervals
between consecutive real singularities. Consider the domain
obtained by slitting $U$ along the two rays $\R\ssm\sigma$ and in
addition delete small disks around singular points as explained in
Lemma~\ref{lem:ind-circ}. The result will be a simply connected
domain $U'$, still containing $T$, whose boundary consists of
three types of arcs/segments (see Fig.~1):
\begin{enumerate}
    \item the ``distant part'', the union of $6$ segments
    belonging to the original boundary of $\partial U$,
    \item ``small arcs'' around singular points $t_j\in\S$
    (semi-circular or circular),
    \item at most $2m$ real segments $\sigma_j\subset\R\ssm\S$
    (upper and lower ``shores'' of the slits).
\end{enumerate}

As before, the number of poles of $f$ in $U'$ is bounded by a
computable function (all of these must be poles of the functions
$x_i\in\C(X)$). Thus to majorize the number of zeros of $f$ in
$T$, it is sufficient to majorize the index of the boundary
$\partial U'$. The contribution from the first two parts (the part
sufficiently distant from the singular locus and small arcs around
Fuchsian singularities) is computable by Lemmas~\ref{lem:ind-seg}
and~\ref{lem:ind-circ}. What remains is to estimate the
contribution of the index of $f$ along the real segments
$\sigma_i$ with endpoints at the singular locus.

7. The index $V_\sigma(f)$ along each such segment
$\sigma=\sigma_j$ by Proposition~\ref{prop:petrov} is bounded in
terms of the number of zeros of $\Im_{\sigma} f$. Computation of
this imaginary part yields (recall that $x_\nu=1$)
\begin{equation*}
  \Im_{\sigma}f=\Im_\sigma h_\nu+\sum_{i=1}^{\nu-1}\Im_\sigma(x_i
  h_i)=0+\sum_{i=1}^{\nu-1}h_i\Im_\sigma
  x_i=\sum_{i=1}^{\nu-1}h_i\~x_i,
\end{equation*}
where the functions $\~x_i=\Im_\sigma x_i$ again belong to $\C(X)$
by Proposition~\ref{prop:reim}.

8. This computation allows to conclude that the contribution of
each segment $\sigma_j$ to the total index of the boundary is
majorized by the number of zeros on $\sigma_j\subset U\ssm\S$ of
an auxiliary function $f_j\in\C(X)[H]$ which is representable as a
sum of $\le\nu-1$ terms of the form \eqref{sum}. By the inductive
assumption, this number is bounded by a computable function of the
degree, the parameters and the distance between $U$ and $\S$.

Thus the total number of zeros of $f$ in $U'$ is shown to be
bounded by a primitive recursive function of admissible parameters
(the algorithm for computing this function is summarized in
\secref{sec:algorithm} below). The proof of Lemma~\ref{lem:isomon}
is complete.\qed

\subsection{Quasialgebraicity near infinity}
The singular point at $t=\infty$ is non-Fuchsian for systems of
the special class, and hence it deserves a separate treatment.

Let $U$ be a \emph{complement to the rectangle},
\begin{equation*}
  U=\C\ssm\{|\Re t|<2,\ |\Im t|<1\}
\end{equation*}
and $\C(Y)$ a Picard--Vessiot field of special type with the given
parameters $n,m,r$ (we retain the notation compatible with that
used earlier in the proof of Lemma~\ref{lem:isomon}). Denote by
$\gamma$ the (interior) boundary of $U$. By definition of the
special class, all finite ramification points of $\C(Y)$ are
outside $U$.

\begin{Lem}\label{lem:isomon-infty}
The field $\C(Y)$ is quasialgebraic in $U$ provided that the
spectral condition \eqref{spec} holds for the ``large loop''
$\gamma$ going around infinity.
\end{Lem}

\begin{proof}
Despite its different appearance, this is again an isomonodromic
reduction principle: the role of the second system is played by an
appropriate Euler system.

1. If the monodromy operator $M$ for the large loop has all
eigenvalues on the unit circle, then there exists an Euler system
\begin{equation*}
  \dot X=t^{-1}A X,\qquad A\in\Mat_n(\C),
\end{equation*}
isomonodromic with the initial system in $U$, whose matrix residue
$A$ has only real eigenvalues in the interval $[0,1)$. Since
conjugacy (gauge transformation)  by constant nondegenerate
matrices does not affect the isomonodromy, one can assume that the
norm of $A$ is bounded, say, by $n$.

2. The Picard--Vessiot field $\C(X)$ is quasialgebraic by
Lemma~\ref{lem:euler-root} globally, hence in $U$ as well.

3. To prove quasialgebraicity of $\C(Y)\subseteq\C(X,H)$, where
$H=YX^{-1}$, we follow the same strategy. More precisely, in this
case choose $U'$ to be the domain $U$ slit along the negative real
semiaxis and with a ``small neighborhood of infinity'' (i.e.,
complement to a disk of large radius) deleted.

As before, we have to majorize the number of roots of functions
from $\C(X)[H]$ via the index of the boundary. This bound will
serve all triangles $T\subset U$ not intersecting the negative
semiaxis (if $T$ intersects the negative semiaxis, the slit must
be along the positive semiaxis).

4. The boundary of $U'$ also consists of three types of
segments/arcs:
\begin{enumerate}
    \item the ``distant part'', the union of $5$ segments
    belonging to the original boundary of $\gamma=\partial U$,
    \item one circular ``small arc'' $\gamma_\infty$ around infinity,
    \item $2$ real segments $\sigma_j\subset\R\ssm\S$
    (upper and lower ``shores'' of the slit).
\end{enumerate}
The only difference from the previous case is that the index of
$\gamma_\infty$ is not guaranteed to be computable, since
Lemma~\ref{lem:ind-circ} does not apply to the infinite singular
point. Yet the following argument shows that for functions from
$\C(X)[H]$ the index is bounded \emph{from above}. In case the
monodromy at infinity is trivial, this would correspond to
bounding the order of pole; note that the order of multiplicity of
a root (the upper bound for the above index) would be much harder
to majorize (cf.~with \cite{bolibr:multiplicity}).

\begin{Prop}\label{prop:index-infty}
If the orientation on $\gamma_\infty$ is chosen counterclockwise
\(as a part of $\partial U'$\), then for any function
$f\in\C(X)[H]\subset\C(X,Y)$ the index $V_{\gamma_\infty}(f)$ is
bounded from above by a computable function of $d=\deg_{X,H} f$
and the parameter $r$ of the special class only.
\end{Prop}

\begin{proof}[Proof of the Proposition]
We claim that any function $f\in\C(X)[H]\subseteq\C(X,Y)$ admits
an asymptotic representation $f(t)=t^\lambda\ln^k (c+o(1))$ with a
real \emph{growth exponent} $\lambda$, natural $k<n$ and a nonzero
complex constant $c\in\C$ near infinity (more precisely, in any
slit neighborhood of infinity). Moreover, we claim that the growth
exponent $\lambda$ is bounded from above in terms of $\deg_{X,H}
f$ and the parameters of the field $\C(X,Y)$.

Indeed, for functions $\C[X]$ this is obvious, since $\C[X]$ as a
ring is generated over $\C$ by monomials $t^\lambda\ln_k t$ with
$\lambda\in[0,1)$ and $k\le n$. For the field $\C(X)$ it follows
from the fact that a ratio of two functions with the specified
asymptotics is again a function with the same asymptotical
representation and the exponent being the difference of the growth
exponents. Thus the growth exponent for $f\in\C(X)$ with $\deg
f=d$ is between $-d$ and $d$.

By definition of the special class, all entries of $Y$ and
$Y^{-1}$ grow no faster than $t^r$ as $t\to\infty$ along any ray.
This means that $H,H^{-1}$ may have a pole of order $\le r+1$ at
infinity, hence the growth exponent of any function $f\in\C[H]$ of
degree $d$ can be at most $(r+1)d$.

Since the growth exponent of a sum is the maximum of the growth
exponents of the separate terms, the growth exponent of any
$f\in\C(X)[H]$ is bounded in terms of $\deg f$ and $r$ as
asserted.

It remains only to note that the index of any function
$f=t^\l\ln^k t(c+o(1)$ along a circular arc with sufficiently
large radius and the angular length less than $2\pi$, is no
greater than $2\pi\l+1$. This proves the Proposition.
\end{proof}

5. The proof of Lemma~\ref{lem:isomon-infty} is completed exactly
as the proof of Lemma~\ref{lem:isomon}. Namely, to majorize
inductively the index of any sum \eqref{sum} along the boundary
$\partial U'$, we estimate the contribution of the ``distant
part'' by Lemma~\ref{lem:ind-seg}, the contribution of the arc
$\gamma$ by Proposition~\ref{prop:index-infty}. The contribution
of two remaining real intervals is majorized by the number of
zeros of two auxiliary functions (the imaginary parts of
$f/x_\nu$) having the same structure \eqref{sum} but with fewer
terms, $\nu-1$.
\end{proof}

\subsection{Algorithm}\label{sec:algorithm}
The above proof of Lemma~\ref{lem:isomon} is algorithmic. Let
$f\in\C(Y)$ be an arbitrary function of degree $d$. The procedure
suggested to majorize the number of its zeros in $U$ (i.e., in any
triangle $T\subset U\ssm\S$) looks as follows.

1$^\circ$. Using the identity $Y=H^{-1}X$, we can express $f$ as
an element of $\C(X,H)$ and then as a ratio of two
``polynomials'', $f=f_1/f_2$, $f_1,f_2\in\C(X)[H]$. The degree
$d'$ of $f_i$ in $X,H$ is immediately computable, as well as the
number of terms $\nu\le d'$ in the representation \eqref{sum}.

2$^\circ$. We start from the list of two functions $\mathcal
L=\{f_1,f_2\}$, both written in the form of a sum \eqref{sum}.
Then we transform this list, repeatedly using the following two
rules:\begin{enumerate}
    \item if $f=\sum_{i=1}^k x_i h_i$ is in $\mathcal L$ and
    $x_k\not\equiv 1$, then $f$ is replaced by $2$ functions, $x_k$ and
    $f/x_k=(\sum_{i=1}^{k-1}x_i/x_k)+h_k$;
    \item if $x_k\equiv 1$, then $f\in\mathcal L$ is replaced by $2m-2$
    functions $\Im_{\sigma_j}f$, the imaginary parts of $f$ on
    each of the real segments $\sigma_j$, $j=1,\dots,2(m-1)$.
\end{enumerate}
These steps are repeated in a simple loop no more than $2d'$
times, after which $\mathcal L$ contains no more than $(2m)^{d'}$
functions, all of them belonging to $\C(X)$.

3$^\circ$. After each cycle, the maximal degree $\max\deg x_i$ of
all coefficients of all functions in $\mathcal L$ is at most
doubled (the imaginary part does not change it, division by the
last term increases at most by the factor of $2$). Thus all
functions obtained after the loop is finished, will have
explicitly bounded degrees not exceeding $d''=d'\cdot2^{d'}$.

4$^\circ$. Since the field $\C(X)$ is quasialgebraic, the number
of isolated zeros of all functions is bounded in terms of their
degrees and the parameters of the field $\C(X)$.

It remains to add the computed bounds together, adding  the
computable contribution coming from the index of all ``distant''
and ``small'' arcs occurring in the process. This contribution is
bounded by a computable function of $d''$, the maximal degree ever
to occur in the construction, the inverse distance from $\partial
U$ to the combined singular locus, and the parameters of the
compositum field $\C(X,Y)$. The algorithm terminates. It will
later be inserted as a ``subroutine'' into one more simple loop
algorithm to obtain the computable bound in
Theorem~\ref{thm:main}.

\section{Isomonodromic surgery}\label{sec:surgery}

Existence of a system isomonodromic to a given one in a given
domain may be problematic in the Fuchsian class, but is always
possible within the special class of systems.

\subsection{Matrix factorization problem}
Let $R=\{a<|t|<b\}\subset\C$ be a circular annulus centered at the
origin on the complex plain and $H(t)$ a matrix function,
holomorphic and holomorphically invertible in the annulus. It is
known, see
\cite{birkhoff-factor,krein-gohberg,markus-feldman,bolibr:kniga},
that $H$ can be factorized as follows,
\begin{equation*}
  H(t)=F(t)\,t^D G(t),\qquad D=\diag\{d_1,\dots,d_n\},
\end{equation*}
where the matrix functions $F$ and $G$ are holomorphic and
holomorphically invertible in the sets $U=\{|t|<b\}\subset\C$ and
$V=\{|t|>a\}\cup\{\infty\}\subset\C P^1$ respectively (their
intersection is the annulus $R$), and $D$ is the diagonal matrix
with the integer entries uniquely defined by $H$.

We need a quantitative version of this result. Unfortunately, it
is impossible to estimate the sup-norms of the matrix factors $F$
and $G$ and their inverses in terms of the matrix norm of $H$, as
we would like to have. A weaker assertion in this spirit still can
be proved.

Assume that a given matrix function $H$, holomorphic and
holomorphically invertible in the annulus $R$, admits
factorization by two matrix functions $F,G$ satisfy the following
properties,
\begin{equation}\label{decomp}
\begin{gathered}
  H(t)=F(t)G(t),
  \\
  \begin{aligned}
  F,F^{-1}\quad&\text{holomorphic and invertible in
  $U=\{|t|<b\}$},
  \\
  G,G^{-1}\quad&\text{holomorphic and invertible in
  $V=\{a<|t|<+\infty\}$}
  \\
  &\text{and have a finite order pole at $t=\infty$}.
  \end{aligned}
\end{gathered}
\end{equation}

A  matrix function $G(t)$, having a unique pole of finite order
$\le \nu$ at infinity, can be represented as follows,
\begin{equation}\label{expand-G}
  G(t)=\~G(t)+G_1 t+G_2 t^2+\cdots+G_\nu t^\nu,\qquad
  G_i\in\Mat_n(\C),
\end{equation}
with $\~G(t)$ holomorphic and \emph{bounded} in $\{|t|>a\}$.
Similarly, if the inverse matrix $G^{-1}(t)$ has a unique pole at
infinity, then
\begin{equation}\label{expand-G'}
  G^{-1}(t)=\~G'(t)+G_1't+G_2't^2+\cdots+G_\nu't^\nu,\qquad
  G_i'\in\Mat_n(\C),
\end{equation}
with a holomorphic bounded term $\~G'(t)$.

We will say that the decomposition \eqref{decomp},
\eqref{expand-G}--\eqref{expand-G'} is \emph{constrained} in a
smaller annulus $R'=\{a'<|t|<b'\}\Subset R$, $a<a'<b'<b$,  by a
finite positive constant $C$, if
\begin{equation}\label{bounds}
\begin{gathered}
  \|F(t)\|+\|F^{-1}(t)\|\le C,\qquad |t|<b,
  \\
  \|\~G(t)\|+\|\~G'(t)\|\le C,\qquad |t|>a,
  \\
  \|G_1\|+\cdots+\|G_\nu\|+\|G_1'\|+\cdots+\|G_\nu'\|+\nu\le C.
\end{gathered}
\end{equation}

Let $R'=\{a'<|t|<b'\}=U'\cap V'\Subset R$ be the annulus,
$a'=a+\tfrac18(b-a)>a$, $b'=b-\tfrac18(b-a)<b$, and
$U'=\{|t|<b'\}$, $V'=\{|t|>a'\}$.

\begin{Lem}\label{lem:kgb}
Assume that the conformal width $b/a-1$ and the exterior diameter
of the annulus $R$ are constrained by the inequalities
\begin{equation}\label{width-parameter}
  b<q,\quad b/a>1+1/q,\qquad q\in\mathbb N.
\end{equation}
Assume that the matrix function $H$ is bounded together with its
inverse in $R$,
\begin{equation}\label{H-bound}
  \|H(t)\|+\|H^{-1}(t)\|\le q',\qquad a<|t|<b,\quad q'\in\mathbb N.
\end{equation}

Then there exist two matrix functions $F$ and $G$, holomorphic and
holomorphically invertible in smaller circular domains $U'$ and
$V'$ respectively, such that the decomposition \eqref{decomp} is
constrained in $R'$ by a bound $C$ given by a computable
\(primitive recursive\) function $\mathfrak C(q,q')$ of the
natural parameters $q,q'$.

If $H$ is real on intersection with the real axis $R\cap\R$, then
$F$ and $G$ can be also found real on $\R$.
\end{Lem}

The proof is given in Appendix~\ref{sec:factorization}.

\subsection{Isomonodromic surgery}
Consider a linear system \eqref{ls}, \eqref{special} from the
special class $\mathcal S(n,m,r)$ and assume that its singular
locus $\S$ consists of two subsets $\S_1,\S_2$ sufficiently well
apart. We prove that the system is isomonodromic to a ``simpler''
system (having fewer singular points), in appropriate simply
neighborhoods of each subset $\S_i$. The proof is constructive and
we show that the new systems can be also chosen from the special
classes with the respective parameters $(n,m_i,r_i)$, $i=1,2$,
explicitly bounded in terms of the parameters of the initial
system and the size of the gap between $\S_1$ and $\S_2$.

\begin{Lem}\label{lem:surgery}
Assume that the singular locus $\S$ of a system \eqref{ls},
\eqref{special} from the special class is disjoint with the
annulus $R=\{a<|t|<b\}$.

Then there exists a system from the special class, isomonodromic
with the initial system \eqref{ls} in the disk $\{|t|<b\}$ and
having no other finite singularities.

The parameters of the new system are bounded by a computable
function of $q$ and the parameters of the initial system.

The monodromy group of the new system satisfies the spectral
condition \eqref{spec} if this condition was satisfied by the
initial system \eqref{ls}.
\end{Lem}

\begin{proof}
The construction is pretty standard, see
\cite{bolibr:kniga,ai:ode-enc}. We have only to verify the
computability of all the relevant constraints.

Let $\gamma$ be a loop running counterclockwise along the middle
circle $\{|t|=\tfrac12(a+b)\}$ of the annulus and $X(t)$ a
solution defined by the condition $X(t_*)=E$ at some positive real
point $t_*$ inside the annulus.

Since the matrix $A(t)$ is explicitly bounded on $\gamma$, the
Gronwall inequality implies that the monodromy factor $M$ of this
solution and its inverse $M^{-1}$ both have the norm explicitly
bounded in terms of the available data, hence the spectrum of $M$
belongs to some annulus $R^*\subset\C$ on the $\l$-plane (actually
we will need only the case when this spectrum belongs to the unit
circle).

The annulus $R^*$ can be slit along a meridian $\{\Arg t=\const\}$
so that the boundary $\gamma^*$ of the resulting simply connected
domain will be sufficiently distant from the spectrum of $M$.
Consider the matrix logarithm defined by the integral
representation, see \cite{lancaster},
\begin{equation*}
  A_0=\frac1{2\pi i}\oint_{\gamma^*}(\l E-M)^{-1}\,\ln\l\,d\l
  \in\Mat_{n\times n}(\C).
\end{equation*}
By construction, $\|A_0\|$ is bounded by a computable function of
the admissible parameters, and the multivalued function
$X'(t)=t^{A_0}$ has the same monodromy along the loop $\gamma$ as
the initial system. Therefore the matrix ratio
$H(t)=X(t)(X'(t))^{-1}$ is single-valued, and holomorphic. By the
Gronwall inequality, $H,H^{-1}$ are bounded in the smaller
sub-annulus $R''=\{a''<|t|<b''\}$, $a''=a+\tfrac14(b-a)$,
$b''=b-\tfrac14(b-a)$, together with the inverse $H^{-1}$ by a
computable function.

Consider the factorization \eqref{decomp} of $H$, described in
Lemma~\ref{lem:kgb} (with $R$ replaced by $R''$). Then the two
multivalued functions $Y(t)=F^{-1}(t)X(t)$ and $Z(t)=G(t)X'(t)$,
the first defined in $U'$ and the other in $V'$, coincide on the
intersection $R'=U'\cap V'$, so that the meromorphic matrix-valued
1-forms $dY\cdot Y^{-1}$ and $dZ\cdot Z^{-1}$, defined in the
domains $\{|t|<b\}$ and $\{|t|>a\}$, in fact coincide on the
intersection $R'$. Hence they are restrictions of some globally
defined meromorphic matrix 1-form $B(t)\,dt$ with a rational
matrix function $B(t)$.

We claim that $B(t)$ has a form \eqref{special} and compute (i.e.,
majorize) its number of singularities $m$ and height $r$. Indeed,
by construction $B(t)$ is real on $R\cap\R_+$ hence everywhere on
$\R$, and
\begin{equation*}
  B(t)=\begin{cases}
  -F^{-1}(t)\cdot\dot F(t)+F^{-1}A(t)F(t),\qquad &|t|<b,
  \\
  \dot G(t)G^{-1}(t)+t^{-1}G(t)A_0G^{-1}(t),\qquad &|t|>a.
  \end{cases}
\end{equation*}
The terms $F^{-1}(t)\cdot\dot F(t)$ and $\dot G(t)\cdot G^{-1}(t)$
are holomorphic at all finite points of their respective domains
of definition. Thus the matrix function $B(t)$ has finite poles
only at those poles of $A(t)$, which are inside the disk
$\{|t|<a\}$. The residue of $B(t)$ at a singular point $t_j$ is
$F^{-1}(t_j)A_jF(t_j)$, where $A_j$ is the residue of $A(t)$ at
this point. Since residues of the initial system were bounded by
the parameters of the special class while the decomposition
\eqref{decomp} was explicitly constrained, all residues of $B(t)$
at all finite points are real and explicitly bounded.

It is easy to verify that the order of pole of $B(t)$ at
$t=\infty$ and its Laurent coefficients are bounded: this follows
from computability of the order of pole of $G$ and $G^{-1}$ at
infinity and the explicit bounds \eqref{bounds} on the matrix
Laurent coefficients $G_i,G_i'$.

It remains to verify the regularity of the singular point at
infinity and compute the growth exponent of the solution $V(t)$.
This is also obvious, since the growth exponents of $G$, $G^{-1}$
are explicitly bounded whereas the growth exponents of $X'(t)$ and
its inverse are bounded by $\|A_0\|$.

The last assertion of the Lemma (on the monodromy group) is
obvious: any simple loop avoiding the part of the spectrum inside
the annulus, can be deformed to remain also inside the annulus,
hence avoiding also singularities outside. The gauge
transformation replacing $X$ by $F^{-1}X$, does not affect the
monodromy factors for the loops on which the transformation is
defined. Hence the spectrum of all such loops remains the same for
the initial and the constructed system.
\end{proof}

\section{Demonstration of the main Theorem~\ref{thm:main}}

In this section we prove our main result, Theorem~\ref{thm:main},
by induction in the number of finite singular points. The idea is
to break the singular locus $\S$ into two parts sufficiently apart
to belong to two disjoint rectangles $U_1,U_2$, and find for each
rectangle a system from the special class (with fewer
singularities), isomonodromic to the given one. Application of the
isomonodromic reduction principle would yield then a computable
global upper bound for the number of zeros of solutions.

\subsection{Stretch and break}\label{sec:stretch}
Given a system \eqref{ls} from the special class \eqref{special}
with $m$ finite singular points, we cover the complex plane $\C$
by three polygonal domains $U_1,U_2,U_3$ (two rectangles and a
complement to a rectangle, all with sides parallel to the
imaginary and real axes) so that in each of these domains the
problem of counting zeros and poles of functions from $\C(X)$ is
reduced to the same problem restated for the three auxiliary
fields $\C(Y_i)$, $i=1,2,3$, each having no more than $m-1$
ramification points.

The first step is to stretch the independent variable so that the
singular locus $\S$ is not collapsing to a single point. By
definition of the standard class, $\S\subset[-1,1]$, and one can
always make an affine transformation (stretch) $t=at'+b$,
$a,b\in\R$, $0<a\le 1$, $|b|\le 1$, which would send the extremal
points $t_{\text{min}}<t_{\text{max}}$ of $\S$ to $\pm 1$.

Such transformation preserves the form \eqref{special} of the
coefficients matrix, eventually affecting only the magnitude of
the Laurent coefficients $A_k'$ at infinity. However, because of
the inequalities $|a|<1$, $|b|<1$ the impact will be limited.
Indeed, the linear transform $t=at'$ with $|a|<1$ can only
decrease the norms $\|A_i'\|$, while the bounded translation
$t=t'+b$, $|b|<1$, may increase them by a factor at most $2^p$.

Thus without loss of generality one can assume that the endpoints
of $\S$ coincide with $\pm 1$. Since there are $m$ singular
points, some of the intervals $(t_j,t_{j+1})$ should be at least
$2/(m-1)$-long, hence the straight line $\ell=\{\Re
t=\tfrac12(t_j+t_{j+1})\}$ subdivides the singular locus into two
parts, each of them at least $1/(m-1)$-distant from $\ell$.

Let $U=\{|\Re t|\le 2,\ |\Im t|<1\}$ be the rectangle containing
inside it the singular locus. Denote by $U_1$ and $U_2$ two parts,
on which $U$ is subdivided by $\ell$, and let $U_3$ be the
complement to $U=U_1\cup U_2$.

For each $\S_i=\S\cap U_i$, $i=1,2$, one can construct the annulus
satisfying the assumptions of Lemma~\ref{lem:surgery} with $q\le
m$ (eventually, translated by no more than 1 along the real axis).
Thus there exist two systems, isomonodromic with the initial
system \eqref{ls}, in the domains $U_1,U_2$, whose boundaries are
at least $1/(m-1)$-distant from all singular points. In $U_3$ the
initial system has only one singularity (at infinity) and hence is
isomonodromic with an appropriate Euler system generating the
corresponding Picard--Vessiot field $\C(Y_3)$.

Notice that all three auxiliary systems have fewer (at most $m-1$)
finite singular points, and all belong to the appropriate special
classes whose parameters are expressed by computable functions of
the parameters of the initial system.

\subsection{Algorithm for counting zeros}
Given a Picard--Vessiot field $\C(X)$ of special type described by
the given set of parameters $n,m,r$ and a function $f\in\C(X)$ of
degree $d$ from this field, we proceed with the following process.

1$^\circ$. Stretch the singular locus $\S$ and subdivide $\C$ into
three domains $U_1,U_2,U_3$ as explained in \secref{sec:stretch},
and for each subdomain construct the auxiliary field $\C(Y_i)$,
$i=1,\dots,3$, isomonodromic to $\C(X)$ in $U_i$. Each is again a
special Picard--Vessiot field with no more than $m-1$ points and
all other parameters explicitly bounded.

2$^\circ$. By the isomonodromic reduction principle
(Lemma~\ref{lem:isomon} for $U_1,U_2$ and
Lemma~\ref{lem:isomon-infty} for $U_3$)  the question on the
number of zeros of $f$ in $U_i$ can be reduced to that for
finitely many functions $f_{i,\alpha}$ from $\C(Y_i)$, each of
them having the degree (relative to the corresponding field) no
more than $2^d$. This step of the algorithm is realized by the
``subroutine'' described in \secref{sec:algorithm} which itself is
an iterated loop. Yet the number of iterations is no greater than
$O(2^d)$.

3$^\circ$. The steps 1$^\circ$ and 2$^\circ$ have to be repeated
for each of the functions $f_{i,\alpha}$ of known degrees in the
respective Picard--Vessiot fields of special type with explicitly
known parameters (actually, only the height has to be controlled).

4$^\circ$. After $m$ loops of the process explained before, all
functions would necessarily belong to the Euler fields. Their
number, their degrees and the heights of the corresponding Euler
fields will be bounded by a computable function of the parameters
of the process. Since the spectral condition \eqref{spec} for each
of these fields will be satisfied, Lemma~\ref{lem:euler-root}
gives an upper bound for the number of isolated zeros of each of
these functions in any triangle free from singularities.

5$^\circ$. Assembling together all these inequalities for the
number of zeros together with the computable contribution from the
index of ``distant'' boundary segments and ``small'' arcs, gives
an upper bound $\mathfrak N'(d,n,m,r)$ for the number of zeros of
the initial function $f$ of degree $d$.

6$^\circ$. The running time of this algorithm is obviously bounded
by an explicit expression depending only on $d,n,m,r$. As is
well-known \cite{meyer-ritchie,machtey}, the bound itself must be
therefore a primitive recursive function of its arguments. This
completes the proof of Theorem~\ref{thm:main}.

\subsection{Fields with pairwise distant ramification
 points}\label{sec:distant}
The arguments described above prove also the following result
described earlier in \secref{sec:distant-intro}.

Assume that the natural number $s\in\N$ is so large that the
singular points $t_1,\dots,t_m$ of the Fuchsian system \eqref{ls},
\eqref{A} satisfy the inequalities
\begin{equation}\label{spread-points}
  |t_i-t_j|\le 1/s,\quad |t_i|\le s,\qquad  s\in\N.
\end{equation}

\begin{Thm}\label{thm:distant}
A Picard--Vessiot field $\C(X)$ of Fuchsian type is
quasialgebraic, provided that its singular points $t_1,\dots,t_m$
satisfy the inequality \eqref{spread-points} and all residues
$A_1,\dots,A_m,A_\infty=-\sum_1^m A_j$ have only real eigenvalues.
Besides the parameters $n,m,r$, the bound for the number of zeros
in this case depends on the additional parameter $s$.
\end{Thm}

\begin{proof}
To prove Theorem~\ref{thm:distant}, it is sufficient to subdivide
$\C$ into several rectangles $U_i$, $i=1,\dots,m$ each containing
only one singularity, with boundaries sufficiently distant from
the singular locus (the parameter $s$ controls this distance from
below). The spectral condition for the small loops $\gamma_j$
around each singular point $t_j$, finite or nor, is automatically
verified: the eigenvalues $\mu_{k,j}$ of the corresponding
monodromy operators $M_j\in\GL(n,\C)$ are exponentials of the
eigenvalues $\l_{k,j}$ of the residues $A_j$, $\mu_{k,j}=\exp 2\pi
i \l_{k,j}$ \cite{bolibr:kniga,ai:ode-enc}.

Application of the isomonodromic reduction principle allows to
derive quasialgebraicity of $\C(X)$ from that for each Euler
field. Note that in this case there is no need in the
isomonodromic surgery and induction in the number of singular
points.
\end{proof}

Since any three points on the sphere $\C P^1$ can always be placed
by an appropriate conformal isomorphism to $0,1,\infty$ so that
\eqref{spread-points} is satisfied with $s=1$, we have the
following corollary.

\begin{Cor}
Any Picard--Vessiot field with $3$ ramification points, is
quasialgebraic, provided that all three residue matrices have only
real eigenvalues.
\end{Cor}

\appendix

\section{Folding of Fuchsian systems}\label{sec:folding}

We first show how to complete the field $\C(\psi^*X)$ to a
Picard--Vessiot field $\C(Y)$ in the case when $\psi(z)=\sqrt z$
is the algebraic map inverse to the \emph{standard fold}
$\f_0\:t\mapsto t^2$. Our main concern is how to place an upper
bound on the height of the ``folded'' system.

Then we construct the rational map $\f$ as an alternating
composition of simple folds and conformal isomorphisms of $\C P^1$
so that it would align all singularities of any given Fuchsian
system \eqref{ls}--\eqref{A} on the real axis, while increasing in
a controllable way the height $r$ and the dimension $n$. On the
final step we symmetrize the system to achieve the condition that
all residues are real.

\subsection{Simple fold}\label{sec:fold}
Consider the Fuchsian system \eqref{ls}--\eqref{A}, assuming that
both $t=0$ and $t=\infty$ (the critical values of the standard
fold $\f_0$) are nonsingular for it:
\begin{equation}\label{0inf}
  t_j\ne 0,\quad j=1,\dots,m,\qquad \sum_{j=1}^m A_j=0.
\end{equation}

We claim that the collection of $2n$ functions
$\bigl(x_1(t),\dots,x_n(t),tx_1(t),\dots,tx_n(t)\bigr)$,
 after the change of the independent
variable $t=\sqrt z$ satisfies a Fuchsian system of $2n$ linear
ordinary differential equations which can be obtained by
separating the even and odd parts $A_\pm(\cdot)$ of the matrix
$A(t)=A_+(t^2)+t\,A_-(t^2)$.

More precisely, the (column) vector function $y(z)=\bigl(x(\sqrt
z),\sqrt z\,x(\sqrt z)\bigr)$ satisfies the system of linear
ordinary differential equations with a rational matrix of
coefficients $B(z)$,
\begin{equation}\label{suspended-system}
\begin{gathered}
    \frac{dy}{dz}=B(z)y,\qquad
    B(z)=\frac12
  \begin{pmatrix}
    A_{-}(z) & z^{-1}A_{+}(z) \\
    A_{+}(z) & A_-(z)+z^{-1}E
  \end{pmatrix},
  \\
  A_+(z)=\sum_{j=1}^m \frac {t_jA_j}{z-z_j},\quad
  A_-(z)=\sum_{j=1}^m\frac{A_j}{z-z_j},\qquad z_j=t_j^2.
\end{gathered}
\end{equation}
This can be verified by direct computation. The explicit formulas
\eqref{suspended-system} allow to compute residues of $B$, showing
that it is indeed Fuchsian under the assumptions \eqref{0inf}.
Note that even if there are two points $t_j=-t_i$ with $z_i=z_j$,
the folded system nevertheless has only simple poles.

The points $z_j=t_j^2$ are finite poles for $B(z)$ with the
corresponding residues $B_j$ given by the formulas
\begin{equation}\label{suspended-residues1}
  B_j=\frac12
  \begin{pmatrix}
    A_{j} & z_j^{-1}A_{j} \\
    z_jA_{j} & A_{j} \
  \end{pmatrix}, \qquad j=1,\dots,m.
\end{equation}
In addition  to these singular points (obviously, simple poles),
the matrix function $B(z)$ has two additional singular points
$z=0$ and $z=\infty$. The residues at these points are
\begin{equation}\label{suspended-residues2}
 B_0=\frac12
  \begin{pmatrix}
    0 & -\sum z_j^{-1}A_j \\
    0 & E \
  \end{pmatrix},
  \qquad
  B_\infty=\frac12\begin{pmatrix}
    0 & 0 \\
    -\sum z_j A_j & -E \\
  \end{pmatrix}.
\end{equation}

Inspection of these formulas immediately proves the following
result.

\begin{Prop}\label{prop:ht-fold}
The transform $\C(\psi^* X)$ of the Picard--Vessiot field of
Fuchsian type satisfying the condition \eqref{0inf} is a subfield
of the Picard--Vessiot field $\C(Y)$ built from solutions of the
system \eqref{suspended-system}.

This system has the double dimension, its singular locus on $\C
P^1$ consists of the squares $z_j=t_j^2$, $j=1,\dots,m$, and the
points $0$ and $\infty$.

The height $r'$ of the folded system \eqref{suspended-system} is
bounded by a computable function of the height $r$ of the initial
system and the inverse distance
$s=\lceil\dist^{-1}(\S,\{0,\infty\})\rceil
=\lceil\max_{j=1,\dots,m}(|t_j|,\,|t_j^{-1}|)\rceil$ between the
polar locus of \eqref{ls} and the critical locus of the fold.
\end{Prop}

Of course, this ``computable function'' can be immediately written
down, but our policy is to avoid explicit formulas.

\subsection{Monodromy of the folded system}
We show now that the spectral condition \eqref{spec} is preserved
by the standard fold. This is not completely trivial.

The proof begins by constructing a fundamental matrix solution
$Y(z)$ to the system \eqref{suspended-system}. Let $z=b$ be a
nonsingular point for the latter, in particular, $b\ne 0,\infty$,
and denote by $t_+(z)$, $t_-(z)$, two branches of the root $\sqrt
z$ for $z$ near $b$, with $t_+(b)=-t_-(b)$ denoted by $a$. Let
$X_\pm(t)$, be two germs of two fundamental matrix solutions for
the system \eqref{ls} defined near the points $\pm a$ (they must
not necessarily be analytic continuations of each other). Then the
formal computations of the previous sections show that (the germ
at $a$ of) the analytic matrix function
\begin{equation}\label{matrixY}
  Y(z)=
  \begin{pmatrix}
    X_+(t_+(z)) & X_-(t_-(z)) \\
    t_+(z)\, X_+(t_+(z)) & t_-(z)\,X_-(t_-(z))
  \end{pmatrix},
  \quad t_\pm\:(\C,b)\to(\C,\pm a),
\end{equation}
satisfies the system \eqref{suspended-system}. In order to show
that it is nondegenerate, subtract from its $(n+k)$th row the
$k$th row multiplied by $t_+(z)$. The result will be a block upper
triangular matrix with two $n\times n$-blocks $X_+(t_+(z))$ and
$(t_-(z)-t_+(z))\,X_-(t_-(z))$ on the diagonal, each of which is
nondegenerate due to the choice of $b\ne0,\infty$.

\begin{Prop}\label{prop:folded-monodromy}
If all simple monodromies of the system \eqref{ls} have
eigenvalues only on the unit circle, then this also holds true for
the folded system \eqref{suspended-system}.
\end{Prop}

\begin{proof}
Consider any simple loop $\gamma$ on the $z$-plane, starting at
the point $a$ and avoiding the folded locus
$\S'=\{0,z_1,\dots,z_m,\infty\}$. Denote by $\gamma_\pm$ its two
preimages on the $t$-plane, obtained by analytic continuation of
the branches $t_\pm(z)$ along $\gamma$. Two cases are to be
distinguished.

If $\gamma$ was not separating $0$ and $\infty$, then both
branches are single-valued on $\gamma$ and hence both $\gamma_\pm$
will be closed loops on the $t$-plane. Denote by $M_\pm$ the
respective monodromy matrices for $X_\pm$ along $\gamma_\pm$.
Clearly, the monodromy of the whole matrix $Y(z)$ along such a
loop will be block-diagonal: $\Delta_{\gamma_\pm}X_\pm=X_\pm
M_\pm$, $\Delta_{\gamma_\pm}t_\pm=t_\pm$, and hence
\begin{equation*}
  \Delta_\gamma
  \begin{pmatrix} X_+&X_-\\t_+ X_+&t_-X_-\end{pmatrix}
    =
  \begin{pmatrix} X_+&X_-\\t_+ X_+&t_-X_-\end{pmatrix}
  \cdot
  \begin{pmatrix}
  M_+&\\&M_-
  \end{pmatrix}.
\end{equation*}
The spectrum of such matrix is the union of spectra of $M_\pm$,
hence the assertion of the Lemma trivially holds.

In the other case when $\gamma$ separates $0$ and $\infty$, the
branches of $t_\pm$ will change sign after going around it, hence
the images $\gamma_\pm=t_\pm(\gamma)$ will be non-closed arcs
connecting two symmetric points $\pm a$. However, making the
second turn around $\gamma$ will restore the closedness: both
$\gamma_+'=\gamma_-\circ\gamma_+$ and
$\gamma_-'=\gamma_+\circ\gamma_-$ will be \emph{simple} closed
loops issued from the points $t_+(a)$ and $t_-(a)$ respectively.
Denote by $M_\pm'$ the monodromy matrices of the solutions $X_\pm$
along the corresponding loops $\gamma_\pm'$. Then, since
$\Delta_{\gamma_\pm'}t_\pm=(-1)^2t_\pm=t_\pm$, the \emph{square}
$\Delta^2_\gamma Y$ of the monodromy factor for the solution $Y$
will be block-diagonal, $\diag(M_+',M_-')$ and hence its spectrum
is on the unit circle.

But if the eigenvalues of the square of a matrix are on the unit
circle, then the eigenvalues of the matrix itself are also of unit
modulus, and the proof of the Lemma is complete.
\end{proof}

\subsection{Placement of poles by a conformal isomorphism}
Applying the simple fold that preserves the real axis can decrease
the number of non-real poles, if some of them were on the
imaginary axis. Thus combining simple folds with real shifts
$t\mapsto t+c$, $c\in \R$ sufficiently many times, one may hope to
place all poles of the resulting system on the real axis.

However, it can happen that after a shift one of the poles would
appear arbitrarily closely to $t=0$ or $t=\infty$ (or simply
coincide with them). This will result in an arbitrarily large
height of the folded system \eqref{suspended-system}.

In order to correct this, instead of shifts (affine translations)
between the folds we will make conformal isomorphisms that will
keep all singularities away from the critical locus $\{0,\infty\}$
of the standard fold. Such isomorphisms can be always found.

Let $\rho=\rho_m=\pi/(2m-2)$ be the constant depending only on the
number $m$ of the points.

\begin{Prop}\label{prop:goodfold}
Let $t_1,\dots,t_m$ be any $m$ points on $\CP^1$ \(real or not\),
and $t_1\notin\R P^1$. Then there exists a conformal isomorphism
of the sphere $\CP^1$ which takes $t_1$ into $i=+\sqrt{-1}$,
preserves the real line and such that the images of the points
$t_2,\dots,t_m$ are all in the annulus $\rho_m\le |z|\le
\rho_m^{-1}$.
\end{Prop}

\begin{proof}
After an appropriate affine transformation preserving $\R$, we may
already assume that $t_1=+i$ while other points $t_j$,
$j=2,\dots,m$ can be arbitrary.

Using the (inverse) stereographic projection, we can identify $\C
P^1$ with the Euclidean sphere $S^2\subset \R^3$ so that the real
line becomes its equator $E$, whereas the pair of points $\pm i$
is mapped into the north and south poles $N,S$ respectively. The
points $0$ and $\infty$ become a pair of antipodal points on $E$.

The conformal isomorphism we are looking for, will be constructed
as a rigid rotation of the sphere around the axis $NS$ so that
some two opposite points of the equator, sufficiently distant from
all other points $t_2,\dots,t_m$, will occupy the positions for
$0$ and $\infty$ respectively. Such a pair of opposite points
always exists for the following simple reason.

Consider $2m-2$ points, $t_2,\dots,t_m$ and their antipodes on the
sphere, and a spherical cap of geodesic radius less than
$\rho=\rho_m=\pi/(2m-2)$ around each point. The union of these
caps cannot completely cover the equator $E$ whose length is
$2\pi$. Therefore there is a point $b$ on the equator, at least
$\rho$-distant from all $a_j$ and their antipodes. Hence both $b$
and anti-$b$ are $\rho$-distant from all $t_2,\dots,t_m$.

It remains to observe that after the rotation sending $b$ to $0$
and anti-$b$ to $\infty$ and returning to the initial affine chart
on $\CP^1$ by the (direct) stereographic projection, the distance
from $0$ to any of $t_j$, equal to $|t_j|$, will be no smaller
than its spherical counterpart, and the same for the distance
$|1/t_j|$ from $t_j$ to infinity. Since the stereographic
projection of the sphere on the plane is conformal (and obviously
so is the rigid rotation of the sphere), the result will be a
conformal isomorphism as required.
\end{proof}

\subsection{Proof of Lemma~\ref{lem:folding}}
We construct an alternate sequence of conformal isomorphisms and
simple folds as follows. Assuming that $k\ge 0$ singularities
already belong to the real axis, choose one that is not in $\R$.
Using the conformal transformation described in
Proposition~\ref{prop:goodfold}, place the chosen point to
$+\sqrt{-1}$ while preserving the real line and keeping all other
singularities away from $\{0,\infty\}$. Then after the standard
fold the number of poles of the system off the real line will be
by one less than before and the construction continues by
induction in $k$.

After at most $m-3$ folds all singular points of the initial
system will be aligned along the real axis (without loss of
generality the first three points can be assumed already real).

The chain of folds and conformal isomorphisms determines the chain
of Fuchsian systems (defined on different copies of $\CP^1$): we
start from the initial system and then either make a conformal
change of the independent variable (the procedure that does not
change neither the residues nor the monodromy group) or replace
the system by the folded system \eqref{suspended-system} of the
double dimension. If the intermediate conformal isomorphisms are
chosen as in Proposition~\ref{prop:goodfold}, then norms of the
residues of the suspended system after the subsequent simple folds
could exceed those of the initial system at most by a constant
factor depending only on $m$. All this follows from
Proposition~\ref{prop:ht-fold}.

The spectral condition \eqref{spec} is also preserved in this
suspension (Proposition~\ref{prop:folded-monodromy}). Thus the
system obtained at the end satisfies all assertions of the
Lemma~\ref{lem:folding} except that its residue matrices are
possibly not real yet. We show that this can be corrected by one
more doubling of dimension, this time involving reflection in the
real axis rather than reflection in the origin.

Together with the system $\dot x=A(t)x$, whose matrix is a
rational function not necessarily real for $t\in\R$, consider its
``mirror image'' $\dot x^*=A^*(t)x^*$ obtained by reflection in
the real axis, denoting by $A^*(t)=\overline{A(\bar t)}$ the
\emph{rational} matrix function and $x^*$ the new dependent
complex vector variable.

Notice that $A^*(t)$ and $A(t)$ take conjugate values for all
$t\in\R$. The variables $u=\tfrac 12 (x+x^*)$ and $v=\tfrac
1{2i}(x-x^*)$, $u,v\in\C^n$, together satisfy the system of
equations
\begin{equation*}
  \begin{pmatrix}
  \dot u\\ \dot v
  \end{pmatrix}
  =\frac12\,\begin{pmatrix}
  A+A^*&i(A-A^*)
  \\
  \noalign{\vskip2pt}
  -i(A-A^*)&A+A^*
  \end{pmatrix}
  \begin{pmatrix}u\\ v\end{pmatrix}
\end{equation*}
and the matrix of this system is real on $\R$.

Applying this last dimension doubling to the system obtained on
the preceding step, we construct a new system with the matrix of
coefficients real on $\R$. This completes the proof of
Lemma~\ref{lem:folding}. \qed

\section{Zeros of functions defined by systems of polynomial
ordinary differential equations}\label{sec:meandering}

In this section we derive Lemmas~\ref{lem:ind-seg}
and~\ref{lem:ind-circ} from the general theorems on zeros of
functions defined by systems of polynomial differential equations
\cite{annalif-99,fields}.

\subsection{A complex analog of the Vall\'ee-Poussin theorem}
 Consider a linear $n$th order differential equation
with analytic coefficients $a_k(t)$ bounded in a domain
$U\subset\C$,
\begin{equation}\label{lode}
\begin{gathered}
  y^{(n)}+a_1(t)\,y^{(n-1)}+\cdots+a_{n-1}(t)\dot y+a_n(t)y=0,
  \\
  |a_k(t)|\le C\qquad \forall t\in U\subset\C,\qquad  k=1,\dots,n.
  \end{gathered}
\end{equation}

\begin{OtherThm}[\cite{fields}, cf.~with \cite{poussin}]
 \label{thm:sturm}
If $y=f(t)\not\equiv0$ is a nonzero solution of the equation
\eqref{lode} and $\gamma=[t_0,t_1]\subset U$ a rectilinear segment
of length $|\gamma|$, then the index $V_\gamma(f)$ of $f$ along
$\gamma$ is explicitly bounded from two sides,
\begin{equation}\label{bd}
  |V_\gamma(f)|\le \pi(n+1)(1+3|\gamma|C).\qed
\end{equation}
\end{OtherThm}

\subsection{Index of small circular arcs around a Fuchsian singularity}
If \eqref{lode} has a singularity (pole of the coefficients) in
$U$, the coefficients become unbounded and Theorem~\ref{thm:sturm}
is not applicable. Yet if the singularity is Fuchsian (i.e., if
the coefficient $a_k$ in \eqref{lode} has a pole of order $\le k$
at this point, see \cite{ai:ode-enc,ince}), then the index can be
majorized along small circular arcs centered at the singular
point.

Let $\mathcal E=t\tfrac d{dt}$ be the Euler differential operator
and $\~a_1(t),\dots,\~a_n(t)$ holomorphic germs at $(\C,0)$.
Consider the differential equation
\begin{equation}\label{eulode}
  \mathcal E^n y+\~a_1(t)\mathcal E^{n-1}y+\cdots+\~a_{n-1}\mathcal Ey+\~a_n(t)y=0,
  \qquad \mathcal E=t\tfrac d{dt},\ t\in(\C,0).
\end{equation}
After expanding the powers $\mathcal E^k$ and division by $t^n$
the equation \eqref{eulode} becomes an equation of the form
\eqref{lode} with coefficients $a_k(t)$ having a pole of order
$\le k$ (i.e., Fuchsian) at the origin $t=0$. Conversely, any
linear equation \eqref{lode} having a Fuchsian singular point at
the origin, can be reduced to the form \eqref{eulode} with
appropriate analytic coefficients $\~a_k(t)$ obtained as linear
combinations of holomorphic germs $t^k\,a_k(t)$ with bounded
constant coefficients.

In the ``logarithmic time'' $\tau=\ln t$ the Euler operator
$\mathcal E$ becomes the usual derivative $\frac d{d\tau}$ and
hence the whole equation \eqref{eulode} takes the form
\eqref{lode} with the coefficients $a_i(\tau)=\~a_i(e^\tau)$
defined, analytic and $2\pi i$-periodic, hence bounded, in a
shifted left half-plane $\{\Re\tau<-B\}$ for some $B\in\R$.
Circular arcs  of angular length $\f$, $0<\f<2\pi$, with center at
the origin, become rectilinear segments of (Euclidean) length
$\f\le 2\pi$ in the chart $\tau$. Application of
Theorem~\ref{thm:sturm} to this equation proves the following
result.

\begin{Cor}\label{cor:smallarc-euler}
Assume that the coefficients of the Fuchsian equation
\eqref{eulode} are explicitly bounded at the origin,
$|\~a_i(0)|\le C$, $i=1,\dots,n$.

Then index $V_\gamma(f)$ of any solution $f\not\equiv0$ along any
sufficiently small circular arc $\gamma\subset\{|t|=\e\}$ of
angular length less than $2\pi$ is explicitly bounded from two
sides,
\begin{equation}\label{bd-arc}
  |V_\gamma(f)|\le \pi(n+1)(1+6\pi C).\qed
\end{equation}
\end{Cor}

A globalization of this result looks as follows. Let
$\S=\{t_1,\dots,t_m\}\subset\C$ be a finite point set, $t_i\ne
t_j$, $\delta$ a monic polynomial with simple roots on $\S$ and
$\mathcal D$ the differential operator of the following form,
\begin{equation}\label{delta}
  \mathcal D=\delta(t)\tfrac d{dt},
  \qquad\delta(t)=(t-t_1)\cdots(t-t_m)\in\C[t].
\end{equation}
Consider the linear differential equation with Fuchsian
singularities on $\S$,
\begin{equation}\label{pseuler}
  \mathcal D^n y+a_1(t)\,\mathcal
  D^{n-1}y+\cdots+a_{n-1}(t)\,\mathcal Dy+a_n(t)y=0,
\end{equation}
where $a_k\in\C[t]$ are polynomials in $t$. We will refer to the
sum of absolute values of (complex) coefficients of a polynomial
$p\in\C[t]$ as its \emph{height}.

\begin{Prop}\label{prop:rat-equat}
If the degrees and heights of $\delta\in\C[t]$ and all polynomial
coefficients $a_k\in\C[t]$ of the equation \eqref{pseuler} do not
exceed $r\in\mathbb N$, then for any nonzero solution $f$ of this
equation\textup:
\begin{enumerate}
    \item if $\gamma$ is a line segment such that
    $\dist(\gamma,\S\cup\infty)\ge 1/s$ for some $s\in\mathbb N$,
    then $|V_\gamma(f)|$ is bounded by a computable function of
    $n,m,r$ and $s$;
    \item if $t_j$ is sufficiently distant from other singularities,
    i.e., $\dist(t_j,\S\ssm\{t_j\})\ge 1/s$ for some $s\in\mathbb
    N$, and $\gamma_\e\subset\{|t-t_j|=\e\}$ is a simple circular
    arc around $t_j$, then $\limsup_{\e\to 0}|V_{\gamma_\e}(f)|$
    is bounded by a computable function of $n,m,r$ and $s$.
\end{enumerate}
\end{Prop}

\begin{proof}
Both assertions follow from Theorem~\ref{thm:sturm} and
Corollary~\ref{cor:smallarc-euler}.

To prove the first assertion, one has to pass from the
representation \eqref{pseuler} to the expanded form
\begin{equation}\label{expanded}
  y^{(n)}+\frac{\~a_1(t)}{\delta(t)}\,y^{(n-1)}
  +\frac{\~a_2(t)}{\delta^2(t)}\,y^{(n-2)}
  +\cdots
  +\frac{\~a_{n-1}(t)}{\delta^{n-1}(t)}\,\dot y
  +\frac{\~a_n(t)}{\delta^n(t)}\,y=0
\end{equation}
and use the explicit lower bound for $|\delta(t)|$ on $\gamma$. To
prove the second assertion assuming for simplicity that $t_j=0$,
one has to substitute $\mathcal D=\~\delta(t)\cdot\mathcal E$,
where $\~\delta=\delta/t$, and reduce \eqref{pseuler} to the form
\eqref{eulode}. The lower bound on $\~\delta(0)$ coming from the
fact that all other singular points are at least $1/s$-distant
from $t=0$, implies then an upper bound for $\~a_j$.
\end{proof}

\subsection{From systems to high order equations}
 \label{sec:chains}
Consider a system \eqref{ls} of linear ordinary differential
equations with rational coefficients and represent its matrix
$A(t)$ in the following form,
\begin{equation}\label{A-red}
  A(t)=\frac1{\delta(t)}\,P(t),\qquad P=\sum_{k=0}^d P_k t^k,\quad
  \delta(t)=\prod_{j=1}^m(t-t_j),
\end{equation}
where $P(t)$ is a matrix polynomial. The height of $P$ is defined
as the total norm of its matrix coefficients $\sum_k \|P_k\|$.
Note that the denominator $\delta(t)$ is a monic polynomial.

The system \eqref{ls}, \eqref{A-red} can be reduced to one high
order linear equation \eqref{pseuler} or \eqref{expanded} in a
variety of ways. Yet only a few of them allow for an explicit
control of the height of coefficients of \eqref{pseuler} in terms
of the height of $P$ and $\delta$ in \eqref{A-red}. One of such
reductions is outlined below (another approach, suggested in
\cite{alexg:thesis}, is based on completely different ideas and
may eventually lead to better bounds).

Let $p_0(t,x)=c_1 x_1+\cdots+c_n x_n$ be an arbitrary linear form
on $\C^n$ with constant coefficients $c_1,\dots,c_n\in\C$.
Consider the sequence of polynomials
$p_i(t,x)\in\C[t,x_1,\dots,x_n]$, $i=1,2,\dots$, defined by the
recurrent formula
\begin{equation}\label{pi}
  p_{i+1}(t,x)=\delta(t)\pd{p_i(t,x)}{t}
  +\sum_{j,k=1}^n \pd{p_i(t,x)}{x_j}P_{jk}(t)x_k
\end{equation}
where $P_{jk}(t)\in\C[t]$ are entries of the polynomial matrix
function $P(t)$. In other words, $p_{i+1}$ is the Lie derivative
of the polynomial $p_i\in\C[t,x]$ along the polynomial vector
field in $\C^{n+1}$,
\begin{equation}\label{poly-vf}
  \delta(t)\pd{}t+\sum_{j,k=1}P_{jk}(t)x_k\pd{}{x_j}
\end{equation}
orbitally equivalent to the system \eqref{ls}, \eqref{A-red}.
Obviously, all $p_i$ are linear forms in the variables
$x_1,\dots,x_n$.

The chain of polynomial ideals
\begin{equation}\label{chain}
  (p_0)\subseteq (p_0,p_1)\subseteq (p_0,p_1,p_2)\subseteq\cdots\subseteq
  (p_0,\dots,p_k)\subseteq\cdots
\end{equation}
in the Noetherian ring $\C[t,x]$ eventually stabilizes. Hence for
some $\ell<+\infty$ one must necessarily have an inclusion
\begin{equation*}
  p_{\ell}\in(p_0,\dots,p_{\ell-1}),
\end{equation*}
meaning that for an appropriate choice of polynomials
$b_1,\dots,b_\ell\in\C[t,x]$, the identity
$p_{\ell}=\sum_{i=1}^\ell b_ip_{\ell-i}$ must hold. Since all
$p_i$ are linear forms in $x_i$, truncation of this identity to
keep only linear terms implies that
\begin{equation}\label{long-eq}
\begin{gathered}
  p_\ell(t,x)+a_1(t)p_{\ell-1}(t,x)+\cdots+a_\ell(t)p_0(t,x)=0,
  \\
  a_i(t)=-b_i(t,0)\in\C[t],
\end{gathered}
\end{equation}
with \emph{univariate} polynomial coefficients $a_i(t)\in\C[t]$.
Restricting this identity on the solution $x(t)$ of the linear
system \eqref{ls}, \eqref{A-red}, we immediately see that its
$k$th derivative $\tfrac{d^k}{dt^k}y$ is the restriction
$p_k(t,x(t))$ on the same solution.  The conclusion is that the
function $y=p_0(t,x(t))$ satisfies the equation \eqref{pseuler}
with the polynomial coefficients $a_i(t)\in\C[t]$, having Fuchsian
singularities only at the roots of $\delta$.

The polynomial identity \eqref{long-eq} is by no means unique
(this concerns both the length $\ell$ and the choice of the
coefficients $a_i$). However, one can \emph{construct} the
identity \eqref{long-eq} of \emph{computable} length $\ell$ and
with coefficients $a_i\in\C[t]$ of \emph{explicitly bounded
heights}.

Assume that a rational matrix function $A=P/\delta$ of the form
\eqref{A-red} has explicitly bounded degree and height:
\begin{equation}\label{ht-A}
  \deg P,\deg \delta\le m,\qquad
  \text{height of $P$, $\delta\le r$}.
\end{equation}

\begin{OtherThm}[\cite{annalif-99}, especially Appendix B]
 \label{thm:chains}
There exist explicit elementary functions, $\mathfrak l(n,m)$ and
$\mathfrak r(n,m,r)$ with the following properties.

For any linear system \eqref{ls}--\eqref{A-red} constrained by
\eqref{ht-A}, and any linear form $p_0(x)$ with constant
coefficients, one can construct polynomial identity between the
iterated Lie derivatives \eqref{long-eq} so that its length and
the degrees of the coefficients are bounded by $\mathfrak l$ and
their heights by $\mathfrak r$ respectively. \qed
\end{OtherThm}

This is a particular case of a more general theorem proved in
\cite{annalif-99} for Lie derivatives along an arbitrary
polynomial (not necessarily linear) vector field in any dimension.
The functions $\mathfrak l$ and $\mathfrak r$ can be explicitly
computed and their growth rates for large values of $n,m,r$
estimated. These growth rates are very large:
\begin{equation}\label{asymptotics}
\begin{gathered}
  \mathfrak l(n,m)\le n^{m^{O(m^2)}},\qquad
  \mathfrak r(n,m,r)\le (2+r)^{\mathfrak m(n,m)},
  \\
  \mathfrak m(n,m,m)\le \exp\exp\exp\exp(4n\ln m+O(1))\qquad\text{ as }n,m\to\infty.
  \end{gathered}
\end{equation}

The enormity of these bounds is another reason why we never
attempted to write explicitly the bounds whose computability is
asserted by our main results (Theorems~\ref{thm:first},
\ref{thm:second}~and ~\ref{thm:main}). However, the simpler
assertion of Theorem~\ref{thm:distant} can eventually be made
explicit with not-too-excessive bounds based on the results of
Grigoriev \cite{alexg:thesis}.

\subsection{Proof of Lemma~\ref{lem:ind-seg}}
Assertion of this Lemma immediately follows from
Theorem~\ref{thm:chains} and Proposition~\ref{prop:rat-equat} in
view of the construction from the preceding section
\secref{sec:chains}. Indeed, a linear system from the special
class $\mathcal S(n,m,r)$ has the coefficients matrix $A(t)$
which, after reducing to the common denominator, can be written as
a ratio $A=P/\delta$ with explicitly bounded height and degree.

Any linear combination $f$ (a function of degree $1$ in the
respective ring $\C[X]$) satisfies then the linear equation
\eqref{pseuler} and the first assertion of
Proposition~\ref{prop:rat-equat} together with the bounds provided
by Theorem~\ref{thm:chains}, yields an upper bound for the index
of $f$ along any segment $\gamma\subset\C\ssm\S$ sufficiently
distant from the singular locus to admit a lower bound for
$|\delta(t)|$ on it.

To treat polynomials of higher degrees, one may either verify by
inspection that the proof of Theorem~\ref{thm:chains} can be
modified to cover also chains starting from a polynomial
$p_0\in\C[t,x]$ of any degree $d$, the bounds depending on $d$ in
an explicit and computable way. An alternative is to notice that
if $x=(x_1,\dots,x_n)$ is solution of the linear system $\dot
x=Ax$, $A=\{a_{ij}(t)\}_{i,j=1}^n$, then the pairwise products
$x_ix_j$ satisfy the linear system
\begin{equation*}
  \tfrac d{dt} x_{ij}=\sum_{k=1}^n (a_{ik}\,x_kx_j+a_{jk}\,x_ix_k),
  \qquad i,j=1,\dots,n,
\end{equation*}
of dimension $n(n+1)/2$, whose coefficients have explicitly
bounded degree and height. This construction can be iterated as
many times as necessary to treat all monomials of degree $\le d$.
Arranging all these systems in the block diagonal matrix, we
obtain a rational linear system $\dot Y_d=A_d(t) Y_d$, such that
the functions of degree $1$ in $\C[Y_d]$ contain all polynomials
of degree $\le d$ from $\C[X]$. The entries of the matrix $A_{d}$
are selected among the entries of $A=A_1$, hence all other
parameters (except for the dimension) will remain the same, as
well and the number and location of the singularities. This
completely reduces the polynomial case to the linear one treated
first.

The index of a rational function $f\in\C(X)$ of degree $d$ does
not exceed the sum of indices of its numerator and denominator.
This observation completes the proof of the Lemma.\qed

\subsection{Proof of Lemma~\ref{lem:ind-circ}}
The same arguments as used above (with the reference to the second
rather than the first assertion of
Proposition~\ref{prop:rat-equat}), prove also the assertion of
Lemma~\ref{lem:ind-circ} in the particular case when the system is
from the special class and the singular point $t_j$ is finite and
at least $1/s$-distant from all other finite singularities. The
bound for the index in this case depends on the additional natural
parameter $s\in\mathbb N$, not allowed in the formulation of the
Lemma.

If we were dealing with a Fuchsian system, this would pose no
problem: by a suitable M\"obius transformation of the sphere $\C
P^1$ one could move any given point away from the rest of $\S$.
Since conformal changes of $t$ preserve the Fuchsian form and do
not affect the matrix residues, the new system will belong to the
same Fuchsian class $\mathcal F(n,m,r)$.

However, the special class is not invariant by M\"obius
transformations (and even affine changes of $t$ may well affect
the condition $\S\subset[-1,1]$). The arguments necessary involve
``non-special'' systems.

Consider the system \eqref{ls} with the coefficient matrix $A(t)$
of the form \eqref{special}, having a simple pole at the origin
$t_1=0$ and make a non-affine conformal map
\begin{equation}\label{pole21}
  z=\frac t{t-c}\iff t=c\frac z{z-1},\qquad 0\ne c\in\R,
\end{equation}
which preserves the origin and takes infinity to the point $z=1$.

The new coefficient matrix $A(z)=A(t)\,dt/dz$ after such
transformation takes the form
\begin{equation}\label{pole-at-1}
  A(z)=\frac{A_{\infty}}{z-1}+\sum_{j=1}^m\frac {A_j}{z-z_j}+
  \sum_{k=2}^{r+2}\frac{B_k}{(z-1)^k}.
\end{equation}

Assume that $|c|<1$. Then the total norm
$\|B_2\|+\cdots+\|B_{r+2}\|$ of the Laurent coefficients $B_k$ at
the point $t=1$ is explicitly bounded in terms of $r$ (the common
upper bound for the height of the initial system \eqref{special}
and the order of its pole at infinity).

Indeed, the (scalar) 1-form $t^k\,dt$ is transformed by the change
of variables \eqref{pole21} into the rational 1-form
$-c^{k+1}\big(1+\frac1{z-1}\big)^k\frac{dz}{(z-1)^2}=c^{k+1}\omega_k$,
where $\omega_k$ is a rational 1-form whose Laurent coefficients
are depending only on $k$ and bounded from above. This means that
each matrix-valued 1-form $t^kA_k'(t)\,dt$ is transformed into a
rational matrix 1-form with arbitrarily small matrix residues,
provided that $|c|>0$ is sufficiently small (depending on $k$). In
other words, the residues $B_k$ in \eqref{pole-at-1} can be made
arbitrarily small by a suitable choice of $c$ small enough.

After the transformation \eqref{pole21} all finite nonzero
singular points $t_j$ will occur at $z_j=t_j/(t_j-c)$. Choosing
$|c|$ sufficiently small, one can arrange them arbitrarily close
to the point $z=1$, in particular, inside the annulus $\{\tfrac12<
|z|\le 2\}$. Then the height of the monic denominator
$\delta(z)=(z-1)^{r+2}\prod_1^m (z-z_j)$ will be then explicitly
bounded.

Finally, the residues $A_j$ at all Fuchsian points remain the
same, thus for all sufficiently small values of $|c|$, the
transformation \eqref{pole21} brings the matrix $A(z)$ of the
system into the form \eqref{pole-at-1} with $\sum\|A_j\|\le r$,
$\sum_k\|B_k\|\le 1$, $z_1=0$ and $\tfrac12\le |z_j|\le 2$ for
$j=2,\dots,m$. After reducing it to the common denominator form
\eqref{A-red}, we obtain a system with explicitly bounded height
of $P$ and $\delta$, which has a singular point $z=0$ at least
$1/2$-distant from all other points which are still on the segment
$[-2,2]$.

Now one can safely reduce this system to a the linear equation
\eqref{pseuler} as described in Theorem~\ref{thm:chains} and apply
the second assertion of Proposition~\ref{prop:rat-equat} to obtain
computable two-sided bounds on the index of small circular arcs on
the $z$-plane around the origin. The images of these arcs on the
$t$-plane will be arbitrarily small circular arcs around $t_j$,
though not necessarily centered at this point. The proof of
Lemma~\ref{lem:ind-circ} is complete.\qed

\section{Quantitative factorization of matrix functions \\
in an annulus}\label{sec:factorization}

As often happened before, the proof with quantitative bounds is
essentially obtained by inspection of existing qualitative proofs,
supplying quantitative arguments when necessary.

The standard proof of matrix factorization theorem that we
``quantize'', is that from \cite{bolibr:kniga}. Namely, for the
matrix function $H(t)$ \emph{sufficiently close} to the (constant)
unity matrix $E$, one can prove existence of \emph{holomorphic and
holomorphically invertible} decomposition $H=FG$, including the
point at infinity (i.e., with both $G$ and $G^{-1}$ bounded in the
outer domain). Then, using approximation by rational functions, we
prove that there exist \emph{meromorphic and meromorphically
invertible} decomposition with explicit bounds on the number of
zeros and poles of the determinants $\det F,\det G$. Finally, we
show that all zeros and poles can be forced to ``migrate'' to the
single point at infinity, while retaining control over all
magnitudes.

Throughout this section $R$ is the annulus $\{a<|t|<b\}$ with
$a,b$ constrained by the natural parameter $q$ as in
\eqref{width-parameter}, and $U,V$ denote respectively the
interior and exterior disks, $U=\{|t|<b\}$,
$V=\{|t|>a\}\subset\C$, so that $R=U\cap V$.

\subsection{Holomorphic factorization of near-identity matrix
functions}

\begin{Lem}\label{lem:small-fact}
There exists a computable function $\mathfrak N(q)$ of one integer
argument $q$ with the following property.

For a given annulus $R$ whose dimensions are determined by the
integer parameter  $q\in\mathbb N$ as in \eqref{width-parameter},
and any holomorphic invertible matrix function $H$ in the annulus,
close to identity enough to satisfy the condition
\begin{equation}\label{small-enough}
  \|H(t)-E\|\le 1/\mathfrak N(q),\qquad t\in R,
\end{equation}
there exists the decomposition $H=FG$ with $F$ and $G$ holomorphic
and holomorphically invertible in the respective circular domains
$U$ and $V$. The matrix functions $F,G$ satisfy the constraints
\begin{equation}\label{kgb-small}
  \|F(t)\|+\|F^{-1}(t)\|\le \mathfrak N(q),\qquad
  \|G(t)\|+\|G^{-1}(t)\|\le \mathfrak N(q)
\end{equation}
in their domains.
\end{Lem}

\begin{proof}[Scheme of the proof]
The identity $H=FG$ can be considered as a functional equation to
be solved with respect to $F$ and $G$ for the given matrix
function $H$. Representing the three functions as $H=E+C$,
$F=E+A$, $G=E+B$, this identity can be rewritten as
\begin{equation}\label{fact}
  C=A+B+AB.
\end{equation}
This equation is nonlinear with respect to the pair of matrix
functions $(A,B)$. Its ``linearization'' is obtained by keeping
only ``first order'' terms,
\begin{equation}\label{lin-fact}
  C=A+B.
\end{equation}
Solution of this linearized equation can be immediately given by
the Cauchy integral. The boundary $\gamma$ of the annulus $R$
consists of two circular arcs, the interior arc $\gamma_a$ and the
exterior arc $\gamma_b$ (properly oriented). Writing
\begin{equation*}
  C(t)=\frac1{2\pi i}\oint_\gamma\frac{C(z)\,dz}{z-t}=\frac1{2\pi
  i}\oint_{\gamma_a}\frac{C(z)\,dz}{z-t}+\oint_{\gamma_b}\frac{C(z)\,dz}{z-t}.
\end{equation*}
The first term which we denote by $B(t)$, is holomorphic in
$V=\{|t|>a\}$. The second term is respectively holomorphic in
$U=\{|t|<b\}$ and after being denoted by $A(t)$, constitutes
together with $B$ a holomorphic solution of the linearized problem
\eqref{lin-fact}. This solution considered as an integral operator
$C(\cdot)\mapsto(A(\cdot),B(\cdot))$ is of explicitly bounded
norm. Indeed, if $\|C(t)\|$ is bounded on $\gamma$ by $1$, then
for each $t\in R$ at least one of the distances
$\dist(t,\gamma_a)$ or $\dist(t,\gamma_b)$ will be no smaller than
$1/2p$, where $q$ is the width parameter of the annulus as
described in \eqref{width-parameter}. This means that the norm of
one of the integrals is explicitly bounded in terms of $q$, but
the second integral is then also bounded by virtue of
\eqref{lin-fact}. Note that this boundedness is specific for the
annulus having two disjoint parts $\gamma_a,\gamma_b$ of the
boundary $\gamma$ (normally the Cauchy operators have unbounded
norms).

The nonlinear equation \eqref{fact} can be now solved by the
Newton method. If $H=E+C$ and $\|C\|\le \e$, then factorization of
$H$ can be reduced to factorization of $(E+A)^{-1}H(E+B)^{-1}$: if
the latter is factorized as $F'G'$ as required, the former can
also be factorized as $(E+A)F'\cdot G'(E+B)$. The new difference
$C'=(E+A)^{-1}(E+C)(E+B)^{-1}-E$ is of \emph{quadratic} order of
magnitude, $\|C'\|\le L\e^2$, where $L=L(q)$ is the norm of the
above described integral operator solving the linearized equation.
The process of linearization and solving the linearized equation
can be continued.

The super-exponential convergence of Newton-type iterations
guarantees that if the original difference $\max_{t\in
R}\|H(t)-E\|$ was smaller than, say, $1/2L(q)$, then the process
converges and the ultimate result, constructed as the infinite
product of matrices fast converging to the identity $E$ in the
respective domains $U,V$, will be of matrix norms explicitly
bounded in terms of $q$ only. The accurate estimates can be
extracted from \cite{bolibr:kniga} if necessary.

If $C$ is real on $R\cap\R$, then clearly the solution of the
linearized equation \eqref{lin-fact} is also real on the real
line, and this condition is preserved in the iterations of the
Newton method. This completes the proof of the Lemma.
\end{proof}

\subsection{Meromorphic solutions of the factorization problem}
The following construction is also standard.

\begin{Prop}
Any matrix function $H$ satisfying the assumptions of
Lemma~\ref{lem:kgb}, can be approximated by a rational matrix
function
\begin{equation}\label{rational-appr}
  M(t)=\sum_{k=-d}^d M_k t^k,\qquad M_k\in \Mat_n(\C),
\end{equation}
so that the product $M(t)H^{-1}(t)$ satisfies the condition
\eqref{small-enough}:
\begin{equation*}
  \|M(t)H^{-1}(t)-E\|\le 1/\mathfrak N(q)
\end{equation*}
in the smaller annulus $R'=\{a'<|t|<b'\}\Subset R$,
$a'=a+\tfrac14(b-a)$, $b'=b-\tfrac14(b-a)$.

The degree $d$ and norms of the matrix coefficients $M_k$ will be
bounded by primitive recursive functions of $q$ and $q'$. If $H$
is real on $R\cap\R$, then $M$ will also be real on $\R$.
\end{Prop}

\begin{proof}
Consider the Laurent expansion for $H$, converging uniformly on
$R'\Subset R$ since $H$ is holomorphic in the annulus $R$. The
rate of this convergence can be explicitly estimated in terms of
$q'$ and $q$, see \eqref{width-parameter} and \eqref{H-bound}.
Thus one can immediately estimate the number of terms of the
Laurent expansion that would be sufficient to keep to satisfy the
required accuracy of the approximation so that
\begin{equation*}
  \|H(t)-M(t)\|\le 1/q'\mathfrak N(q).
\end{equation*}
Then $H^{-1}M-E$ and $MH^{-1}-E$ will be as small as asserted on
$R'$. The rest is obvious.
\end{proof}

This Proposition together with Lemma~\ref{lem:small-fact} applied
to $H^{-1}M$, guarantees for an arbitrary holomorphically
invertible matrix function $H$ existence of a factorization of the
form
\begin{equation}\label{merom}
\begin{gathered}
  H=FG, \qquad G=G'M,
  \\
  \text{ $F,G'$ holomorphic in $U',V'$ respectively,
  $M$ rational},
  \end{gathered}
\end{equation}
where $F,G'$ are holomorphic, holomorphically bounded and
invertible in $U'=\{|t|<b'\}\Subset U$, $V'=\{|t|>a'\}\Subset V$,
and are bounded together with their inverses there by a computable
function of $q,q'$ only.

\subsection{Expulsion of poles of $G^{-1}$. Proof of Lemma~\ref{lem:kgb}}
The term $F$ in the decomposition \eqref{merom} already meets the
requirements imposed in \eqref{bounds}.

To achieve the conditions required from the second term $G$, we
multiply it from the left by several rational matrix functions, at
the same time multiplying $F$ from the right by their inverses.
The construction, due to G.~Birkhoff \cite{birkhoff-factor}, is
utterly classical.

The matrix function $G$ has a unique pole of order $\le d$ at
infinity, but may have a number of ``\emph{zeros}'' (poles of the
inverse matrix $G^{-1}$). Their number and position are determined
by zeros of $\det G=\det G'\det M$, that is, by zeros of the
rational function $\det M(t)$ of known degree, since the term $G'$
is holomorphically invertible. Thus $G^{-1}$ has no more than a
computable number of poles $z_1,\dots,z_k$, $k$ being bounded by a
computable function of $q,q'$. We are interested only in those
poles that are in $V'$.

Let $t=z$ be any root of the determinant $\det M(t)$. The matrix
$G(z)$ has nontrivial left null space, since $\det G(z)=0$.
Consider the constant matrix $C=C_z$ built as follows: its first
row is the (left) null vector of $G(z)$ of unit Hermitian norm,
while the other rows are chosen among the rows of the identity
matrix $E$ in such a way that both $\|C\|$ and $\|C^{-1}\|$ are
bounded by $n$. This is always possible: it is sufficient to
delete from the coordinate vectors in $\C^n$ the vector which has
the largest Hermitian scalar product with the null vector.

By construction, the product $C_zG(t)$ at the point $t=z$ has the
zero first row. Hence the product $\diag\{r_z(t),1,\dots,1\}\cdot
CG(t)$, $r_z(t)=|z|/(t-z)$, continues to be holomorphic at the
point $z$. The correction term
$Q_z(t)=\diag\{r_z(t),1,\dots,1\}\cdot C$ is bounded together with
its inverse in the thinner annulus $R''=\{a''<|t|<b''\}\Subset
R'$, the ``middle belt'' of $R'$, $a''=a'+\tfrac14(b'-a')$,
$b''=b'-\tfrac14(b'-a')$. Note that this bound is \emph{uniform}
as $z$ varies in $V'=\{b'<|z|\}$. Besides, $Q_z(t)$ is bounded
together with its inverse also in $U'$. All these bounds are
computable functions of $q,q'$. Replacing the decomposition $H=FG$
by $H=FQ_z^{-1}Q_zG$, we obtain a new factorization of $H$ that
has fewer poles of the right term in $V'$ (both terms remain
holomorphic and the left is holomorphically invertible in $U'$).

This expulsion of a pole of $G^{-1}$ from $t=z$ to $t=\infty$ has
to be repeated for each pole $z_j\in V'$  of the initial matrix
$G=G'M$ (the order is inessential). As a result, we replace the
decomposition \eqref{merom} by
\begin{equation}\label{merom-1}
  H=FQ^{-1}\cdot QG,\qquad \text{$Q=Q(t)$ rational without poles
  or zeros in $R'$}.
\end{equation}
By construction, both $QG$ and $(QG)^{-1}$ have only one
(multiple) pole at $t=\infty$, and all four matrices $FQ^{-1}$,
$QG$ and their inverses, are explicitly bounded in $R''$ in the
sense of the norm by a computable function of $q,q'$.

The order of pole of $QG$ at infinity is no greater than the
initial order $d$ of pole of $M$. The order of pole of $(QG)^{-1}$
does not exceed $k$, the degree of $\det M$ (the number of
singular points that were expelled). Finally, the norms of the
Laurent coefficients of $(QG)^{\pm 1}$ at infinity can be
majorized using the Cauchy estimates (differentiating the Cauchy
integral along the inner boundary of the annulus $R''$). These
bounds will prove the assertion of Lemma~\ref{lem:kgb}.

The only remaining problem is to choose $Q$ being real on $\R$
provided that $G$ is real there. In this is the case, then
non-real ``zeros'' of $G$ come in conjugate pairs $z,\bar z$ and
the corresponding left null vectors (rows) $v,\bar v\in\C^n$ are
complex conjugate. In this case both zeros can be expelled if the
corresponding rational matrix factor $Q_{z,\bar z}(t)$ has the
first row of the form
\begin{equation*}
  |z|\(\frac v{t-z}+\frac{\bar v}{t-\bar z}\)\in\R^n,\qquad t\in\R.
\end{equation*}
This vector function has the norm bounded from above away from the
pair of points $\{z,\bar z\}$ uniformly over all $z\in V'$, and as
before the inverse $Q^{-1}$ is of bounded norm provided that the
other rows of $Q$ were appropriately chosen among the rows of the
identity matrix $E$ exactly as before.

This remark concludes the proof of Lemma~\ref{lem:kgb}.\qed


\def\BbbR{$\mathbf R$}\def\BbbC{$\mathbf
  C$}\providecommand\cprime{$'$}\providecommand\mhy{--}\font\cyr=wncyr9
\providecommand{\bysame}{\leavevmode\hbox
to3em{\hrulefill}\thinspace}
\providecommand{\MR}{\relax\ifhmode\unskip\space\fi MR }
\providecommand{\MRhref}[2]{%
  \href{http://www.ams.org/mathscinet-getitem?mr=#1}{#2}
} \providecommand{\href}[2]{#2}

\end{document}